\documentclass[12pt,letterpaper]{amsart}
\usepackage{amsmath,txfonts}
\usepackage{amssymb}
\usepackage{amsxtra}
\usepackage{amsthm, color}
\usepackage{txfonts}
\usepackage{graphicx}
\usepackage{times}
\usepackage{citeref}
\usepackage{tikz}
\usepackage{hyperref}
\usepackage{stmaryrd}
\usepackage[T3,T1]{fontenc}
\usepackage{pgfplots}

\usetikzlibrary{calc}

\usepackage{mathrsfs}
\usepackage{amsfonts}
\usepackage{amssymb}
\usepackage{ifthen}
\usepackage{graphicx}
\nonstopmode \numberwithin{equation}{section}

\newtheorem{thm}{Theorem}[section]
\newtheorem{lem}{Lemma}[section]%[thm]
\newtheorem{cor}[thm]{Corollary}
\newtheorem{prop}[thm]{Proposition}

\newtheorem{step}{Step}[section]%[thm]

\theoremstyle{definition}
\newtheorem{mlem}{Main lemma}[section]
\newtheorem{assertion}{Assertion}[section]
\newtheorem{cl}{Claim}[section]
\newtheorem{ca}{Case}[section]
\newtheorem{sca}{Subcase}[section]
\newtheorem{scl}{Subclaim}[section]
\newtheorem{conj}[thm]{Conjecture}
\newtheorem{fact}{Fact}[section]
\newtheorem{defn}[thm]{Definition}
\newtheorem{op}[thm]{Open Problem}

\newtheorem{ques}{Question}[section]
\newtheorem{rem}[thm]{Remark}
\newtheorem{exam}[thm]{Example}

\numberwithin{equation}{section}

\newcounter {own}
\def\theown {\thesection       .\arabic{own}}

\newenvironment{pf}[1][]{%
 \vskip 3mm
 \noindent
 \ifthenelse{\equal{#1}{}}%
  {{\slshape Proof. }}%
  {{\slshape #1.} }%
 }%
{\qed\bigskip}

\newcounter{alphabet}
\renewcommand{\thealphabet}{\Alph{alphabet}}

\newenvironment{Thm}[1][]{\refstepcounter{alphabet}%
	\bigskip
	\noindent
	{\bf Theorem \thealphabet}%
	\ifthenelse{\equal{#1}{}}{}{ (#1)}%
	{\bf .} \itshape
}{\vskip 8pt}

\newenvironment{Lem}[1][]{\refstepcounter{alphabet}%
	\bigskip
	\noindent
	{\bf Lemma \thealphabet}%
	{\bf .} \itshape
}{\vskip 8pt}

\def\be{\begin{equation}}
\def\ee{\end{equation}}

\newcommand{\ben}{\begin{enumerate}}
\newcommand{\een}{\end{enumerate}}

\newcommand{\blem}{\begin{lem}}
\newcommand{\elem}{\end{lem}}
\newcommand{\bthm}{\begin{thm}}
\newcommand{\ethm}{\end{thm}}
\newcommand{\bcor}{\begin{cor}}
\newcommand{\ecor}{\end{cor}}
\newcommand{\beg}{\begin{exam}}
\newcommand{\eeg}{\end{exam}}
\newcommand{\begs}{\begin{examples}}
\newcommand{\eegs}{\end{examples}}
\newcommand{\bdefe}{\begin{defn}}
\newcommand{\edefe}{\end{defn}}
\newcommand{\bques}{\begin{ques}}
\newcommand{\eques}{\end{ques}}
\newcommand{\bei}{\begin{itemize}}
\newcommand{\eei}{\end{itemize}}
\newcommand{\bcon}{\begin{conj}}
\newcommand{\econ}{\end{conj}}
\newcommand{\bop}{\begin{op}}
\newcommand{\eop}{\end{op}}

\newcommand{\bas}{\begin{assertion}}
\newcommand{\eas}{\end{assertion}}

\newcommand{\bfa}{\begin{fact}}
\newcommand{\efa}{\end{fact}}

\newcommand{\bca}{\begin{ca}}
\newcommand{\eca}{\end{ca}}

\newcommand{\bst}{\begin{step}}
\newcommand{\est}{\end{step}}

\newcommand{\bsca}{\begin{sca}}
\newcommand{\esca}{\end{sca}}

\newcommand{\bcl}{\begin{cl}}
\newcommand{\ecl}{\end{cl}}

\newcommand{\bmlem}{\begin{mlem}}
\newcommand{\emlem}{\end{mlem}}

\newcommand{\bscl}{\begin{scl}}
\newcommand{\escl}{\end{scl}}

\newcommand{\bcons}{\begin{conjs}}
\newcommand{\econs}{\end{conjs}}

\newcommand{\bprop}{\begin{prop}}
\newcommand{\eprop}{\end{prop}}

\newcommand{\br}{\begin{rem}}
\newcommand{\er}{\end{rem}}
\newcommand{\brs}{\begin{rems}}
\newcommand{\ers}{\end{rems}}
\newcommand{\bo}{\begin{obser}}
\newcommand{\eo}{\end{obser}}
\newcommand{\bos}{\begin{obsers}}
\newcommand{\eos}{\end{obsers}}
\newcommand{\bpf}{\begin{pf}}
\newcommand{\epf}{\end{pf}}
\newcommand{\ba}{\begin{array}}
\newcommand{\ea}{\end{array}}
\newcommand{\beq}{\begin{eqnarray}}
\newcommand{\beqq}{\begin{eqnarray*}}
\newcommand{\eeq}{\end{eqnarray}}
\newcommand{\eeqq}{\end{eqnarray*}}

\newcounter{minutes}
\divide\time by 60
\newcounter{hours}
\multiply\time by 60 \addtocounter{minutes}{-\time}
\begin{document}
	\title{Growth type theorems  of harmonic  $(K,K_{0})$-quasiregular mappings}

%\author{Shaolin Chen${}^{~\mathbf{*}}$}
%\address{S. L. Chen,    Center for Applied Mathematics of Guangxi, Guangxi Normal University,
%Guilin, Guangxi 541004, People's Republic of China} \email{mathechen@126.com}

\author{Shaolin Chen}
 \address{S. L. Chen, Center for Applied Mathematics of Guangxi, Guangxi Normal University,
Guilin, Guangxi 541004, People's Republic of China; College of Mathematics and
	Statistics, Hengyang Normal University, Hengyang, Hunan 421008, People's Republic of China} \email{mathechen@126.com}

%\\$^{\mathbf{*}}$Corresponding author}
	

	\maketitle
	
	%\def\thefootnote{}
	%\footnotetext{
		%\texttt{\tiny File:~\jobname .tex,
			%          printed: \number\year-\number\month-\number\day,
			%          \thehours.\ifnum\theminutes<10{0}\fi\theminutes}
		%}
	\makeatletter\def\thefootnote{\@arabic\c@footnote}\makeatother

\begin{abstract}
The purpose of this paper is twofold.
First, we establish several sharp Hardy-Littlewood type radial growth theorems for harmonic $(K,K_{0})$-quasiregular mappings.
Second, we prove some sharp coefficient growth theorems for such mappings.
In particular, we prove the validity of Das and Kaliraj's conjecture for  certain  univalent harmonic subclasses, namely \( S_{H}(K) \) and \( S_{H}^{0}(K,K_{0}) \).
\end{abstract}

\maketitle \pagestyle{myheadings} \markboth{ S. L. Chen}{Growth type theorems   of harmonic  $(K,K_{0})$-quasiregular mappings}

\section{Introduction}\label{csw-sec1}

For $b\in\mathbb{C}$ and $r>0$, let $\mathbb{D}(b,r)=\{z:~|z-b|<r\}$. In particular, we use $\mathbb{D}_{r}$ to
denote the disk $\mathbb{D}(0,r)$ and $\mathbb{D}$ to denote the unit disk $\mathbb{D}_{1}$. Moreover, let $\mathbb{T}:=\partial\mathbb{D}$ be the
unit circle. We use  ``${\rm Re}$'' and ``${\rm Im}$'' to denote  the real and imaginary parts of a complex number, respectively.
For a continuously differentiable complex-valued function $f$ and a point $z=x+iy\in\mathbb{D}$,
the directional derivative in the direction $\theta\in[0,2\pi]$ is defined as
$$\partial_{\theta}f(z)=\lim_{\rho\rightarrow0^{+}}\frac{f(z+\rho e^{i\theta})-f(z)}{\rho}=f_{z}(z)e^{i\theta}
+f_{\overline{z}}(z)e^{-i\theta},
$$
where $\rho\in(0,1-|z|)$,  $f_{z}:=\partial f/\partial z=1/2\left(\partial f/\partial x-i\partial f/\partial y\right)$ and $f_{\overline{z}}:=\partial f/\partial \overline{z}=1/2\left(\partial f/\partial x+i\partial f/\partial y\right)$.
Then
$$\Lambda_{f}(z):=\max_{\theta\in[0,2\pi]}\big\{|\partial_{\theta}f(z)|\big\}=|f_{z}(z)|+|f_{\overline{z}}(z)|
$$
and
$$\lambda_{f}(z):=\min_{\theta\in[0,2\pi]}\big\{|\partial_{\theta}f(z)|\big\}=\big||f_{z}(z)|-|f_{\overline{z}}(z)|\big|.
$$

\subsection*{$(K,K_{0})$-quasiregular mappings}
A mapping $f:~\Omega\rightarrow\mathbb{C}$ is said to be absolutely
continuous on lines, $ACL$ in brief, in the domain $\Omega$ if for
every closed rectangle $R\subset\Omega$ with sides parallel to the
axes $x$ and $y$, $f$ is absolutely continuous on almost every
horizontal line and almost every vertical line in $R$. Such a
mapping has, of course,  partial derivatives $f_{x}$ and $f_{y}$
a.e. in $\Omega$. Moreover, we say $f\in ACL^{2}$ if $f\in ACL$ and
its partial derivatives are locally $L^{2}$ integrable in $\Omega$.

A   mapping $f$ of  $\Omega$ into
$\mathbb{C}$ is called a {\it $(K,K_{0})$-quasiregular mapping} if

\ben
\item  $f$ is $ACL^{2}$ in  $\Omega$ and $J_{f}\geq0$ a.e. in   $\Omega$, where $J_{f}$ is the  Jacobian of $f$;

\item there are constants $K\geq1$
and $K_{0}\geq0$ such that
$$\Lambda_{f}^{2}\leq KJ_{f}+K_{0}~\mbox{ a.e. in}~\Omega.$$
\een  In particular, if $K_{0}\equiv0$, then a
$(K,K_{0})$-quasiregular mapping is said to be  $K$-quasiregular. Furthermore, we say $f$ is $K$-quasiconformal in $\Omega$ if $f$ is  $K$-quasiregular and homeomorphic in $\Omega$ (see \cite{C-K,C-P-W,C-P,FS,gt,Ni}).
%It is well known that every  quasiregular mapping  is an  quasiregular
%mapping. But the inverse of this statement is not true.

\subsection*{Harmonic mappings and harmonic  $(K,K_{0})$-quasiregular mappings} A two times continuously differentiable complex-valued function $f=u+iv$ in a domain $\Omega\subseteq\mathbb{C}$ is called  harmonic mapping (i.e., a complex-valued harmonic function)
 if the real-valued functions $u$ and $v$ satisfy
 the Laplace equation $\Delta u=\Delta v=0$, where $\Delta$ represents the Laplacian operator
 $$\Delta:=4\frac{\partial^{2}}{\partial z\partial \overline{z}}=
 \frac{\partial^{2}}{\partial x^{2}}+\frac{\partial^{2}}{\partial y^{2}}$$
 and $z=x+iy\in\Omega$ (see \cite{Du}).
 For a function $\psi\in C(\mathbb{T})$, we denote by $P[\psi]$ the Dirichlet solution, for Laplace's operator
$\Delta$, of $\psi$ over $\mathbb{D}$, that is $\Delta P[\psi]=0$ in $\mathbb{D}$ and $P[\psi]=\psi$ on $\mathbb{T}$.
It is well-known that, for  $\psi\in L^{1}(\mathbb{T})$,

%Let us recall that the Poisson integral, $P[\varphi]$, of a function $\varphi\in L^{1}(\mathbb{T})$ is defined by

$$P[\psi](z)=\frac{1}{2\pi}\int_{0}^{2\pi}\psi(e^{{\rm i}\tau})\mathbf{P}(z,e^{{\rm i}\tau})d\tau,$$
where $\mathbf{P}(z,e^{{\rm i}\tau})=\frac{1-|z|^{2}}{|e^{{\rm i}\tau}-z|^{2}}$
is the Poisson kernel. The theory of complex-valued harmonic function has many important applications.
For example, recently, Aleman and Constantin \cite{A-C}, and  Constantin and Mart\'in \cite{C-M}  use it to solve some fluid mechanics problems.

For a sense-preserving harmonic mapping $f$ defined in $\mathbb{D}$,
the Jacobian of $f$ is given by
$$J_{f}=|f_{z}|^{2}-|f_{\overline{z}}|^{2}=\Lambda_{f}\lambda_{f}
$$
and the second complex dilatation of $f$ is given by
$\mu=\overline{f_{\overline{z}}}/f_{z}$.  It is well-known that every harmonic
mapping $f$ defined in a simply connected domain $\Omega$ admits a decomposition $f = h + \overline{g}$, where $h$ and $g$ are analytic;
this decomposition is unique up to an additive constant.
Recall that $f$ is
sense-preserving in $\Omega$ if $J_{f}>0 $ in $\Omega$.
Thus, $f$ is locally univalent and sense-preserving in $\Omega$
if and only if $J_{f}>0$ in $\Omega$; or equivalently if $h'\neq
0$ in $\Omega$ and  $\mu =g'/h'$ has the
property that $|\mu|<1$ in $\Omega$ (see
\cite{Clunie-Small-84,Du,Lewy}). A mapping $f$ in $\Omega$ is called a {\it harmonic $(K,K_{0})$-quasiregular mapping}  if
 $f$ is a harmonic and $(K,K_{0})$-quasiregular mapping in $\Omega$, where $K\geq1$
and $K_{0}\geq0$  are constants. In particular, if $K_{0}=0$, then a harmonic  $(K,K_{0})$-quasiregular mapping in $\Omega$ is called a harmonic $K$-quasiregular
mapping in $\Omega$. Moreover,  we say $f$ is harmonic $K$-quasiconformal in $\Omega$ if $f$ is harmonic  $K$-quasiregular and homeomorphic in $\Omega$.
  %Throughout this paper, we use $\mathscr{H}$ to denote all harmonic mappings of $\mathbb{D}$ into $\mathbb{C}$.
%Furthermore, denote by $\mathscr{A}$ the set of  all analytic functions of $\mathbb{D}$ into $\mathbb{C}$.

\subsection*{The generalized Beltrami
equation} Let $f$ be a  continuously differentiable complex-valued function  in a domain $\Omega\subseteq\mathbb{C}$. The following equation is called {\it the generalized Beltrami
equation} (see \cite{FS,Ni}):
\begin{equation}\label{Bel-q}
\overline{f_{\overline{z}}(z)}=\kappa_{1}(z)f_{z}(z)+\kappa_{2}(z),~z\in\Omega,
\end{equation}
where $\kappa_{1}$ and $\kappa_{2}$ are complex-valued functions in $\Omega$ with
$$\|\kappa_{1}\|_{\infty}=\sup_{z\in\Omega}|\kappa_{1}(z)|<1~\mbox{and}~\|\kappa_{2}\|_{\infty}=\sup_{z\in\Omega}|\kappa_{2}(z)|\geq0.$$
In particular, if $\kappa_{2}=0$ in $\Omega$, then (\ref{Bel-q})
is called {\it the  Beltrami equation} (see \cite{AIM}). Furthermore,
if $\kappa_{1}$   is analytic in $\Omega$, and $\kappa_{2}=0$ in $\Omega$, then (\ref{Bel-q})
is called {\it the second Beltrami equation} (see \cite{Du}). It is easy to see that if $\kappa_{1}=\kappa_{2}=0$ in $\Omega$, then $f_{\overline{z}}=0$ which is {\it the Cauchy-Riemann equation}.

\begin{prop}\label{Prop-1}
Let  $f$ be a  continuously differentiable complex-valued function  in  $\Omega$. If $f$ is a solution to (\ref{Bel-q}) with $J_{f}\geq0$ in $\Omega$,
then $f$ is a $(K,K_{0})$-quasiregular mapping, where $K$ only depends  on $\kappa_{1}$, and $K_{0}$ only depends  on $\kappa_{1}$ and $\kappa_{2}$.
%In particular, if $\kappa_{1}$ and $\kappa_{2}$ are analytic in (\ref{Bel-q}), then a two times continuously differentiable complex-valued function $f$ satisfying (\ref{Bel-q}) is %a harmonic $(K,K_{0})$-quasiregular mapping.
\end{prop}

The main purpose of this paper is to prove sharp theorems on the radial growth (of Hardy-Littlewood type) and on the coefficient growth of harmonic  $(K,K_{0})$-quasiregular mappings.

For convenience, we make a notational convention.
Throughout this paper, we use  $\mathscr{A}$ to denote all  analytic functions of $\mathbb{D}$ into $\mathbb{C}$, and
use $\mathscr{H}$ to denote all harmonic mappings of $\mathbb{D}$ into $\mathbb{C}$. Moreover,
 we use the symbol $C$ to denote various positive
constants, whose values may change from one occurrence to another. Also, we use the notation $C=C(a,b,\ldots)$, which means that the constant $C$ depends only on the given parameters $a$, $b$, $\ldots$.

\section{Preliminaries and  main results}\label{csw-sec1-1}

The {\it  Hardy type space}
$\mathbf{H}^{p}_{g}$ $(p\in(0,\infty])$ consists of all those functions
$f$ of $\mathbb{D}$ into $\mathbb{C}$ such that $f$ is measurable, the integral mean $M_{p}(r,f)$ exists for all $r\in[0,1)$,
$$\|f\|_{p}:=\sup_{r\in[0,1)}M_{p}(r,f)<\infty$$ for
$p\in(0,\infty)$, and $$\|f\|_{\infty}:=\sup_{r\in[0,1)}M_{\infty}(r,f)<\infty$$ for $p=\infty$,
where
$$M_{p}(r,f)=\left(\frac{1}{2\pi}\int_{0}^{2\pi}|f(re^{i\theta})|^{p}d\theta\right)^{\frac{1}{p}}~\mbox{and}~M_{\infty}(r,f)=\sup_{\theta\in[0,2\pi]}|f(re^{i\theta})|.$$
 In particular, we use $\mathbf{h}^{p}=\mathbf{H}^{p}_{g}\cap\mathscr{H}$ and
$H^{p}:=\mathbf{H}^{p}_{g}\cap\mathscr{A}$
to denote the {\it harmonic Hardy space} and the {\it analytic Hardy space}, respectively. If
$f\in\mathbf{h}^{p}$ for some $p\in [1,\infty)$,
then the radial limits
$$f(\zeta)=\lim_{r\rightarrow1^{-}}f(r\zeta)$$ exist for almost every $\zeta\in\mathbb{T}$, and $f\in L^{p}(\mathbb{T})$ (see \cite[Theorems 6.7,  6.13 and 6.39]{ABR}).
  Since $|f|^{p}$ is subharmonic for $p\geq1$, we see that $M_{p}(r,f)$ is increasing on $r\in[0,1)$, and
 \be\label{fx-1}
 \|f\|_{p}^{p}=\lim_{r\rightarrow1^{-}}M_{p}^{p}(r,f)=\frac{1}{2\pi}\int_{0}^{2\pi}|f(e^{i\theta})|^{p}d\theta.
 \ee

\subsection*{Radial growth type theorems   of harmonic  $(K,K_{0})$-quasiregular mappings}
 Let $\chi_{1}(r)$ and $\chi_{2}(r)$ be nonnegative functions on $r \in [0,1)$. Throughout this paper, we write
$$
\chi_{1}(r) = O\bigl(\chi_{2}(r)\bigr)
$$
to denote that $\chi_{1}$ is asymptotically bounded above by $\chi_{2}$ as $r \to 1^{-}$. That is, there exist constants $C > 0$ and $r_{0} \in (0,1)$ such that for every $r \in (r_{0}, 1)$,
$$
\chi_{1}(r) \leq C \chi_{2}(r).
$$

%\subsection*{Conjugate functions type   growth  theorems of harmonic quasiregular mappings}
It is well known that there is a close relation between the integral means of an analytic function and those of its derivative.
For $f\in\mathscr{A}$, a classical result of Hardy and Littlewood asserts that if $p\in(0,\infty]$ and $\beta\in(1,\infty)$, then
$$M_{p}(r,f')=O\left(\frac{1}{(1-r)^{\beta}}\right)$$  if and only if
$$M_{p}(r,f)=O\left(\frac{1}{(1-r)^{\beta-1}}\right)$$  (see \cite[Theorem  5.5]{duren} and \cite{DMP,G-P,H-L,H-L-1}).
In particular, if $\beta=1$, then Hardy and Littlewood established the following result.

\begin{Thm}{\rm (\cite[Theorem  6]{H-L})}\label{Thm-HL-3}
Let $f$ be analytic in $\mathbb{D}$. For $p\in(0,\infty)$, if $$M_{p}(r,f')=O\left(\frac{1}{1-r}\right)$$ 
then for $p\in(0,1]$, $$M_{p}(r,f)=O\left(\left(\log\frac{1}{1-r}\right)^{\frac{1}{p}}\right)$$ 
and for $p\in(1,\infty)$,  $$M_{p}(r,f)=O\left(\log\frac{1}{1-r}\right).$$ 
\end{Thm}

In 2023,  Das and Kaliraj  \cite{DK} extended Hardy-Littlewood's result to harmonic mappings for $\beta\in(1,\infty)$, which is as follows.

\begin{Thm}\label{DK-23}{\rm (\cite[Theorem  1]{DK})}
Let $p\in[1,\infty)$ and $\beta>1$. If $f\in\mathscr{H}$, then
$$M_{p}(r,\Lambda_{f})=O\left(\frac{1}{(1-r)^{\beta}}\right)$$ 
if and only if $$M_{p}(r,f)=O\left(\frac{1}{(1-r)^{\beta-1}}\right).$$ 
\end{Thm}

 A continuous non-decreasing and unbounded function $\eta:~[0,1)\rightarrow(0,\infty)$ is called a
 weight. A weight $\eta$ is called doubling if there is a constant $C>1$ such that
 $$\eta(1-s/2)<C\eta(1-s)$$ for $s\in(0,1]$ (see \cite{AD}).
 Recently, the radial growth properties of integral means of analytic functions have attracted much attention of many authors
(see \cite{DMP,G-P,PP,SS-05}).
In the following,  we extend and improve Theorems \ref{Thm-HL-3} and  \ref{DK-23} into the following form.

\begin{thm}\label{thm-0.3}
Let $f=u+iv$ be  a harmonic $(K,K_{0})$-quasiregular mapping in $\mathbb{D}$, where $K\geq1$ and $K_{0}\geq0$ are constants.
\begin{enumerate}
\item[{\rm $(\mathscr{C}_{1})$}] For $p\in(0,1]$, if $$M_{p}(r,\nabla u)=O\left(\frac{1}{1-r}\right),$$ 
then \be\label{Sh-1}M_{p}(r,f)=O\left(\left(\log\frac{1}{1-r}\right)^{\frac{1}{p}}\right).\ee  In particular, if $p\in[1/2,1)$, then the estimate of (\ref{Sh-1}) is sharp.
\item[{\rm $(\mathscr{C}_{2})$}] For $p\in(1,2)$, if
\be\label{gfk-01}M_{p}(r,\nabla u)=O(\eta(r)),
\ee
then
\be\label{Sh-11}M_{p}(r,f)=O\left(\Psi(r)\right),\ee  where $$\Psi(r)=\left(\int_{0}^{r}(1-\rho^{2})^{p-1}\big(\eta(\rho)\big)^{p}d\rho
+\frac{2^{p-1}}{p}\big(\eta(r)\big)^{p}(1-r)^{p}\right)^{\frac{1}{p}}.$$
On the contrary, if
(\ref{Sh-11}) holds for a doubling function $\Psi$, then
\be\label{eq-101}M_{p}(r,\nabla u)=O\left(\frac{\Psi(r)}{1-r}\right).
\ee

\item[{\rm $(\mathscr{C}_{3})$}] For $p\in[2,\infty)$, if
\be\label{Sh-13}M_{p}(r,\nabla u)=O(\eta(r)),\ee then
\be\label{Sh-14}M_{p}(r,f)=O\left(\Phi(r)\right),\ee where $$\Phi(r)=\left(\int_{0}^{r}\left(1-\rho\right)(\eta(\rho))^{2}d\rho\right)^{\frac{1}{2}}.$$
Conversely, if (\ref{Sh-14}) holds for a doubling function $\Phi$, then
\be\label{eq-102}M_{p}(r,\nabla u)=O\left(\frac{\Phi(r)}{1-r}\right).
\ee 

\end{enumerate}
\end{thm}

\begin{cor}\label{cor-2.2}
If we take $\eta(r)=1/(1-r)$ in (\ref{gfk-01}), then (\ref{Sh-11}) is  refined as follows
\be\label{103-eq}M_{p}(r,f)=O\left(\left(\log\frac{1}{1-r}\right)^{\frac{1}{p}}\right),\ee  and this estimate is sharp.
Moreover,  if, for $q>1$, we take $\eta(r)=1/(1-r)^{q}$ in (\ref{gfk-01}), then (\ref{gfk-01}) holds if and only if
\beqq
M_{p}(r,f)=O\left(\frac{1}{(1-r)^{q-1}}\right).
\eeqq 

Moreover,  if we take $\eta(r)=1/(1-r)$ in (\ref{Sh-13}), then (\ref{Sh-14}) is  refined as follows
\be\label{104-eq}M_{p}(r,f)=O\left(\left(\log\frac{1}{1-r}\right)^{\frac{1}{2}}\right),\ee  and this estimate is sharp.
Furthermore,  if, for $q>1$, we take $\eta(r)=1/(1-r)^{q}$ in (\ref{Sh-13}), then (\ref{Sh-13}) holds if and only if
\beqq
M_{p}(r,f)=O\left(\frac{1}{(1-r)^{q-1}}\right).
\eeqq 
\end{cor}

%A weight $\omega$ is called a {\it Dini weight} if $w(t)/t$ is integrable on $[0,2)$
%A continuous increasing function $\omega:[0,2)\rightarrow[0,\infty)$ with $\omega(0)=0$
%is called a {\it Dini weight} if $\omega(t)/t$ is non-increasing for $t\in(0,2)$ and
%there is a positive constant $C$ such that
%\be\label{eq2x}
%\int_{0}^{\sigma}\frac{\omega(t)}{t}\,dt\leq C \omega(\sigma),~\sigma\in[0,2).
%\ee

Let $E$ be a family of increasing
 concave functions $\varphi:[0,\infty)\rightarrow[0,\infty)$ with $\varphi(0)=0$ (see \cite{Ai}).
 %It follows from \cite[Lemma 2.2]{Ai} that $\omega(t)/t$ is nonincreasing for $t\in(0,\infty)$, where $\omega\in E$.
A concave function $\varphi\in E$ is called a {\it Dini type weight}, if there exist some $\delta_{0}>0$
and a positive constant $C$, such that
\be\label{eq2x}
\int_{0}^{\delta}\frac{\varphi(t)}{t}\,dt\leq C \varphi(\delta)
\ee for $0<\delta<\delta_{0}$ (see   \cite{Dy2}).
For $\varphi\in E$, we use  $\mathscr{L}_{\varphi,p}(\mathbb{D})$ to denote the class of all Borel functions
$f$ of $\mathbb{D}$ into $\mathbb{C}$ such that, for $z_{1},~z_{2}\in\mathbb{D}$,
$$\mathcal{L}_{p}[f](z_{1},z_{2})\leq C\varphi(|z_{1}-z_{2}|),$$
where $C$ is a positive constant and
$$\mathcal{L}_{p}[f](z_{1},z_{2})=
\begin{cases}
\displaystyle\left(\int_{0}^{2\pi}|f(e^{i\eta}z_{1})-f(e^{i\eta}z_{2})|^{p}d\eta\right)^{\frac{1}{p}}
& \mbox{if } p\in(0,\infty),\\
\displaystyle|f(z_{1})-f(z_{2})| &\mbox{if } p=\infty.
\end{cases}
$$
The Lipschitz constant of $f\in\mathscr{L}_{\varphi,p}(\mathbb{D})$ is defined as follows $$\|f\|_{\Lambda_{\varphi,p}(\Omega),s}:=\sup_{z_{1},z_{2}\in\Omega,z_{1}\neq\,z_{2}}\frac{\mathcal{L}_{p}[f](z_{1},z_{2})}{\varphi(|z_{1}-z_{2}|)}<\infty.$$
 Moreover, for $\varphi\in E$, we define the space $\mathscr{L}_{\varphi,p}(\mathbb{T})$ consisting of those $f\in L^{p}(\mathbb{T})$
for which
$$\mathcal{L}_{p}[f](\xi_{1},\xi_{2})\leq C\varphi(|\xi_{1}-\xi_{2}|),~\xi_{1},~\xi_{2}\in\mathbb{T},$$
where $C>0$ is a constant and
$$\mathcal{L}_{p}[f](\xi_{1},\xi_{2})=
\begin{cases}
\displaystyle\left(\int_{0}^{2\pi}|f(e^{i\eta}\xi_{1})-f(e^{i\eta}\xi_{2})|^{p}d\eta\right)^{\frac{1}{p}}
& \mbox{if } p\in(0,\infty),\\
\displaystyle|f(\xi_{1})-f(\xi_{2})| &\mbox{if } p=\infty.
\end{cases}
$$

For  some related studies of $\mathscr{L}_{\varphi,p}(\mathbb{D})$, we refer the readers to \cite{CH-2023-SCM,Dy1,Dy2,P-99,P-08} for details.

Generally speaking, it is reasonable  to expect an analytic function to be smooth on the boundary if its derivative grows slowly, and
conversely. However, Hardy and Littlewood proved the principle which can be expressed in surprisingly precise form, as follows (see \cite[Theorem  5.1]{duren}).

\begin{Thm}{\rm (Hardy-Littlewood's Theorem)}\label{HL-20}
Let $f\in\mathscr{A}$ and let $f$ be continuous in $\overline{\mathbb{D}}$.
Then for
$\alpha\in(0,1]$, we have $f\in\mathscr{L}_{\varphi_{\alpha},\infty}(\mathbb{T})$ if and only if
$$|f'(z)|=O\left(\frac{1}{(1-r)^{1-\alpha}}\right),~z\in\mathbb{D},$$ where $r=|z|$ and $\varphi_{\alpha}(t)=t^{\alpha}$ for $t\geq0$.
\end{Thm}

In \cite{DK}, Das and Kaliraj generalized Theorem \ref{HL-20} to harmonic mappings as follows.

\begin{Thm}\label{Thm-oj}
Let $f\in\mathscr{H}$ and let $f$ be continuous in $\overline{\mathbb{D}}$.
Then  $f\in\mathscr{L}_{\varphi_{\alpha},\infty}(\mathbb{T})$ if and only if
$$\Lambda_{f}(z)=O\left(\frac{1}{(1-r)^{1-\alpha}}\right),
$$  where $r=|z|$ and $\varphi_{\alpha}$ is defined in Theorem \ref{HL-20}.
\end{Thm}

%In \cite[p. 598]{Dy3}, Dyakonov   proved that the following original Hardy-Littlewood theorem  (see also \cite[Ineq. (2.3)]{Dy2}).

%\begin{Thm}\label{CorA-HL}%{\rm (Hardy-Littlewood's Theorem)}
%Suppose that $u$ is a real-valued harmonic function in $\mathbb{D}$ with
% $u\in\mathscr{L}_{\omega_{\alpha},\infty}(\mathbb{D})$, where $\alpha\in(0,1]$ and $\omega_{\alpha}$ is defined in Theorem J.
%  Let $v$ be a harmonic conjugate of $u$ with $v(0)=0$. Then $v\in\mathscr{L}_{\omega_{\alpha},\infty}(\mathbb{D})$.
%\end{Thm}

In the following,
we extend and improve Theorems \ref{HL-20} and \ref{Thm-oj} to harmonic $(K,K_{0})$-quasiregular mappings by employing different proof techniques.

\begin{thm}\label{eq-thm-9}
%Let $\varphi$ be a Dini type weight such that there is a positive constant $C$ such that for all $\delta\in[0,\pi]$,
%$$\delta\int_{\delta}^{\pi}\frac{\varphi(t)}{t^{2}}dt\leq\,C\varphi(\delta).$$ 

Let $\varphi$ be a weight function of Dini type satisfying the following condition: there exists a positive constant $C$ such that for every $\delta \in [0,\pi]$,
$$
\delta \int_{\delta}^{\pi} \frac{\varphi(t)}{t^{2}}\, dt \leq C \varphi(\delta).
$$
Let $f=u+iv$ be  a harmonic $(K,K_{0})$-quasiregular mapping in $\mathbb{D}$ that is continuous in $\overline{\mathbb{D}}$, where $K\geq1$ and $K_{0}\geq0$ are constants.
Then, for $p\in[1,\infty]$,  $u\in\mathscr{L}_{\varphi,p}(\mathbb{T})$ if and only if
\be\label{gv-1}
\begin{cases}
\displaystyle\left(\int_{0}^{2\pi}\left(\Lambda_{f}(ze^{i\theta})\right)^{p}d\theta\right)^{\frac{1}{p}}\leq\,C(K,p)\frac{\varphi\big(1-|z|\big)}{1-|z|}+C(p)\sqrt{K_{0}},
& \mbox{if } p\in[1,\infty),\\
\displaystyle \Lambda_{f}(z)\leq C(K)\frac{\varphi\big(1-|z|\big)}{1-|z|}+C\sqrt{K_{0}} &\mbox{if } p=\infty,
\end{cases}
\ee where $C$, $C(K)$ and $C(K,p)$  are positive constants.

%\begin{enumerate}
%\item[{\rm $(\mathscr{D}_{1})$}] Let $f=u+iv$ be  a harmonic $(K,K_{0})$-quasiregular mapping in $\mathbb{D}$ that is continuous in $\overline{\mathbb{D}}$, where $K\geq1$ and $K_{0}\geq0$ are %constants.
%Then, for $p\in[1,\infty]$,  $u\in\mathscr{L}_{\varphi,p}(\mathbb{T})$ if and only if
%\be\label{gv-1}
%\begin{cases}
%\displaystyle\left(\int_{0}^{2\pi}\left(\Lambda_{f}(ze^{i\theta})\right)^{p}d\theta\right)^{\frac{1}{p}}\leq\,C(K,p)\frac{\varphi\big(1-|z|\big)}{1-|z|}+C(p)\sqrt{K_{0}},
%& \mbox{if } p\in[1,\infty),\\
%\displaystyle \Lambda_{f}(z)\leq C(K)\frac{\varphi\big(1-|z|\big)}{1-|z|}+C\sqrt{K_{0}} &\mbox{if } p=\infty,
%\end{cases}\ee where $C$, $C(K)$ and $C(K,p)$  are positive constants.
%\item[{\rm $(\mathscr{D}_{2})$}] Let $f\in\mathscr{H}$ and let $f$ be continuous in $\overline{\mathbb{D}}$.  Then, for $p\in[1,\infty]$,  $f\in\mathscr{L}_{\varphi,p}(\mathbb{T})$ if and only if %(\ref{gv-1}) holds.
%\end{enumerate}
\end{thm}

%If we replace $u\in\mathscr{L}_{\varphi,p}(\mathbb{T})$ in Theorem \ref{eq-thm-9} by $f\in\mathscr{L}_{\varphi,p}(\mathbb{T})$, then we have the following result.

%\begin{cor}\label{eq-cor-9}
%Suppose that $f\in\mathscr{H}$ and $\varphi$ is a Dini type weight such that there is a positive constant $C$ such that for all $\delta\in[0,\pi]$,
%$$\delta\int_{\delta}^{\pi}\frac{\varphi(t)}{t^{2}}dt\leq\,C\varphi(\delta).$$
 % Then, for $p\in[1,\infty]$, $f$ is continuous in $\overline{\mathbb{D}}$ and $f\in\mathscr{L}_{\varphi,p}(\mathbb{T})$, if and only if there is a positive
 % constant $C$ such that
 % \be\label{gv-1}
%\left(\int_{0}^{2\pi}\left(\Lambda_{f}(ze^{i\theta})\right)^{p}d\theta\right)^{\frac{1}{p}}\leq\,C\frac{\varphi\big(d_{\mathbb{D}}(z)\big)}{d_{\mathbb{D}}(z)},
%\ee where $d_{\mathbb{D}}(z)=1-|z|$.
%\end{cor}

If we take $\varphi(t):=\varphi_{\alpha}(t)=t^{\alpha}~(t\geq0)$ in Theorem \ref{eq-thm-9}, then we obtain the following result, where $\alpha\in(0,1]$.

\begin{cor}
Let $\varphi_{\alpha}(t)=t^{\alpha}$ for $t\geq0$, where  $\alpha\in(0,1)$. Let $f=u+iv$ be  a harmonic $(K,K_{0})$-quasiregular mapping in $\mathbb{D}$ that is continuous in $\overline{\mathbb{D}}$, where $K\geq1$ and $K_{0}\geq0$ are constants.
Then, for $p\in[1,\infty]$, $u\in\mathscr{L}_{\varphi_{\alpha},p}(\mathbb{T})$ if and only if

\be\label{KO-1}
\begin{cases}
\displaystyle\left(\int_{0}^{2\pi}\left(\Lambda_{f}(ze^{i\theta})\right)^{p}d\theta\right)^{\frac{1}{p}}=O\left(\frac{1}{(1-r)^{1-\alpha}}\right),~
& \mbox{if } p\in[1,\infty),\\
\displaystyle \Lambda_{f}(z)=O\left(\frac{1}{(1-r)^{1-\alpha}}\right),~
 &\mbox{if } p=\infty,
\end{cases}
\ee
where $r=|z|$.

%\begin{enumerate}
%\item[{\rm $(\mathscr{E}_{1})$}] Let $f=u+iv$ be  a harmonic $(K,K_{0})$-quasiregular mapping in $\mathbb{D}$ that is continuous in $\overline{\mathbb{D}}$, where $K\geq1$ and $K_{0}\geq0$ are %constants.
%Then, for $p\in[1,\infty]$, $u\in\mathscr{L}_{\varphi_{\alpha},p}(\mathbb{T})$ if and only if

%\be\label{KO-1}
%\begin{cases}
%\displaystyle\left(\int_{0}^{2\pi}\left(\Lambda_{f}(ze^{i\theta})\right)^{p}d\theta\right)^{\frac{1}{p}}=O\left(\frac{1}{(1-|z|)^{1-\alpha}}\right)~
%\mbox{as}~ |z|\rightarrow1^{-},
%& \mbox{if } p\in[1,\infty),\\
%\displaystyle \Lambda_{f}(z)=O\left(\frac{1}{(1-|z|)^{1-\alpha}}\right)~
%\mbox{as}~ |z|\rightarrow1^{-}, &\mbox{if } p=\infty.
%\end{cases}
%\ee

%\item[{\rm $(\mathscr{E}_{2})$}] Let $f\in\mathscr{H}$ and let $f$ be continuous in $\overline{\mathbb{D}}$.  Then, for $p\in[1,\infty]$,  $f\in\mathscr{L}_{\varphi_{\alpha},p}(\mathbb{T})$ if a%nd only if (\ref{KO-1}) holds.
%\end{enumerate}

%Suppose that $f=u+iv$ is  a harmonic $(K,K_{0})$-quasiregular mapping in $\mathbb{D}$, where $K\geq1$ and $K_{0}\geq0$ are constants.
%Then, for  $\alpha\in(0,1]$ and $p\in[1,\infty]$, $u$ is continuous in $\overline{\mathbb{D}}$ and $u\in\mathscr{L}_{\varphi_{\alpha},p}(\mathbb{T})$, if and only if
%$$
%\left(\int_{0}^{2\pi}\left(\Lambda_{f}(ze^{i\theta})\right)^{p}d\theta\right)^{\frac{1}{p}}=O\left(\frac{1}{(1-|z|)^{1-\alpha}}\right)
%$$ as $|z|\rightarrow1^{-}$.
\end{cor}

\subsection*{Coefficient growth type theorems of univalent harmonic  $(K,K_{0})$-quasiregular mappings}
Let ${\mathcal S}_{H}$ denote the family of sense-preserving planar
harmonic univalent mappings $f=h+\overline{g}$ in $\mathbb{D}$
satisfying the normalization $h(0)=g(0)=h'(0)-1=0$, where $h$ and
$g$ are analytic in $\mathbb{D}$.  The family ${\mathcal S}_{H}$
together with a few other geometric subclasses, originally
investigated in detail by \cite{Clunie-Small-84,Small}, became
instrumental in the study of univalent harmonic mappings (see
\cite{Du}) and has attracted much attention of many function
theorists. If the co-analytic part $g$ is identically zero in the
decomposition of $f=h+\overline{g}$, then the class ${\mathcal
S}_{H}$ coincides with the classical family $\mathcal S$ of all
normalized analytic univalent functions
$h(z)=z+\sum_{n=2}^{\infty}a_{n}z^{n}$ in $\mathbb{D}$. If
${\mathcal S}_H^{0}=\{f=h+\overline{g} \in {\mathcal S}_H: \,g'(0)=0
\} $, then the family ${\mathcal S}_H^{0}$ is both normal and
compact (cf. \cite{Clunie-Small-84, Du}).
Denote by ${\mathcal
S}_{H}(K, K_{0})$ (resp. ${\mathcal S}_{H}^{0}(K, K_{0})$) if $f\in {\mathcal
S}_{H}$ (resp. ${\mathcal S}_{H}^{0}$) and is a harmonic $(K,K_{0})$-quasiregular mapping in $\mathbb{D}$,
where $K\geq 1$ and $K_{0}\geq0$ are constants. In particular, if $K_{0}=0$, then we use  ${\mathcal
S}_{H}(K)$ and ${\mathcal S}_{H}^{0}(K)$ to denote ${\mathcal
S}_{H}(K, K_{0})$ and ${\mathcal S}_{H}^{0}(K, K_{0})$, respectively.
  Obviously, ${\mathcal S}={\mathcal S}_{H}^{0}(1)={\mathcal S}_{H}(1)$. In particular, let
$\gamma=\sup_{f\in{\mathcal
S}_{H}}(|f_{zz}(0)|/2)$. It is well known that $$2\leq\gamma\leq\frac{32\pi}{27}(\pi+6\sqrt{3})-2<48.4,$$ but the sharp value of
$\gamma$ is still unknown. For mappings of the class ${\mathcal S}_{H}$, it is conjectured that $\gamma=3$  (see \cite{Clunie-Small-84, Du,Small}).

In \cite{DK}, Das and Kaliraj established a radial growth type estimate and a coefficient bound for $f\in{\mathcal
S}_{H}$ as follows.

\begin{Thm}{\rm (\cite[Theorem 7]{DK})}\label{DK-2023}
Suppose that $p\in[1,\infty)$ and $f=h+\overline{g}\in{\mathcal
S}_{H}$ with $h(z)=z+\sum_{n=2}^{\infty}a_{n}z^{n}$ and $g(z)=\sum_{n=1}^{\infty}b_{n}z^{n}$.
Then for every $\epsilon>0$,
$$M_{p}(r,f)=O\left(\frac{1}{(1-r)^{k(p)+\epsilon}}\right),$$ where $k(p)=\sqrt{\gamma^{2}-\frac{1}{p}+\frac{1}{4p^{2}}}-\frac{1}{2p}$.
Consequently, there is an absolute constant $C$ such that
$$|a_{n}|\leq Cn^{\gamma-\frac{1}{2}}$$ for  $n\in\{2,3,\ldots\}$.
\end{Thm}

Let us recall a famous  generalized Bieberbach conjecture  to the class $ {\mathcal
S}_{H}$. For $f(z)=z+\sum_{n=2}^{\infty}a_{n}z^{n}+\sum_{n=1}^{\infty}b_{n}z^{n}\in{\mathcal
S}_{H},$ it is conjectured that

$$|a_{n}|<\frac{1}{3}\left(2n^{2}+1\right)~\mbox{and}~|b_{n}|<\frac{1}{3}\left(2n^{2}+1\right),$$
which implies that $\gamma=3$ and
$$\sup_{n\geq2}\left\{\frac{|a_{n}|}{n^{\gamma-1}}\right\}<\infty.$$ This conjecture has been hanging in the air until now
(see \cite{Du}). Inspired by this conjecture, Das and Kaliraj in \cite{DK} posed the following conjecture.

\begin{conj}{\rm (\cite[Concluding remarks]{DK}})
Under the hypothesis of Theorem E, there is an
absolute constant $C$ such that $$|a_{n}|\leq Cn^{\gamma-1}$$ for  $n\in\{2,3,\ldots\}$.
\end{conj}

%In the following, we prove this conjecture is true for  $f\in{\mathcal S}_{H}(K)$ and  $f\in{\mathcal S}_{H}^{0}(K, K_{0})$, respectively.

In the following, we investigate the validity of Das and Kaliraj's conjecture for  certain  univalent harmonic subclasses, namely $ S_{H}(K) $ and $ S_{H}^{0}(K,K_{0}) $.

\begin{thm}\label{thm-01}
Suppose that  $f=h+\overline{g}\in{\mathcal
S}_{H}(K)$ with $h(z)=z+\sum_{n=2}^{\infty}a_{n}z^{n}$ and $g(z)=\sum_{n=1}^{\infty}b_{n}z^{n}$, where $K\geq1$ is a constant.
Let $\gamma^{\ast}=\sup_{f\in{\mathcal
S}_{H}(K)}(|f_{zz}(0)|/2)$.
 Then, for $r\in(0,1)$ and $p\in(1/\gamma^{\ast},\infty)$,
\be\label{eq-th-c1}M_{p}(r,f)=O\left(\frac{1}{(1-r)^{\gamma^{\ast}-\frac{1}{p}}}\right),\ee
and there exists a constant  $C(K)$ such that
\be\label{eq-th-2}|a_{n}|\leq C(K)n^{\gamma^{\ast}-1}\ee for  $n\in\{2,3,\ldots\}$.
In particular, when $K=1$,  the exponent $\gamma^{\ast}-1$ in (\ref{eq-th-2}) reduces to $1$.
\end{thm}

The following result easily follows from Theorem \ref{thm-01}.

\begin{cor}
Suppose that  $f=h+\overline{g}\in{\mathcal
S}_{H}(K)$ with $h(z)=z+\sum_{n=2}^{\infty}a_{n}z^{n}$ and $g(z)=\sum_{n=1}^{\infty}b_{n}z^{n}$, where $K\geq1$ is a constant.
 Then, for $r\in(0,1)$ and $p\in(1/\gamma,\infty)$,
\beqq\label{eq-th-cx1}M_{p}(r,f)=O\left(\frac{1}{(1-r)^{\gamma-\frac{1}{p}}}\right),\eeqq
and there exists a constant $C(K)$ such that
\be\label{eq-th-x2}|a_{n}|\leq C(K)n^{\gamma-1}\ee for  $n\in\{2,3,\ldots\}$.
In particular, when $K=1$, the exponent $\gamma-1$ in (\ref{eq-th-x2}) reduces to $1$.
\end{cor}

\begin{thm}\label{thm-02}
Suppose that
$f=h+\overline{g}\in{\mathcal S}_{H}^{0}(K, K_{0})$ with $h(z)=z+\sum_{n=2}^{\infty}a_{n}z^{n}$ and $g(z)=\sum_{n=2}^{\infty}b_{n}z^{n}$, where
$K\geq1$ and $K_{0}\geq0$  are constants. Then, for $r\in(0,1)$ and $p\in(1/\gamma,\infty)$,
\be\label{eq-th-t1}M_{p}(r,f)=O\left(\frac{1}{(1-r)^{\gamma-\frac{1}{p}}}\right),\ee
and there exists a constant $C(K,K_{0})$ such that
\be\label{eq-th-t2}|a_{n}|\leq C(K,K_{0})n^{\gamma-1}\ee for  $n\in\{2,3,\ldots\}$.
In particular, when $K-1=K_{0}=0$, the exponent $\gamma-1$ in (\ref{eq-th-t2}) reduces to $1$.
\end{thm}

%The proofs of , Theorems   \ref{maint}, \ref{thm-0.1}, \ref{thm-K-L} and \ref{thm-07} will be presented in Sec. \ref{csw-sec2}.
%The proof of Proposition \ref{Prop-1}, Theorems \ref{thm-0.3}, \ref{eq-thm-9}, \ref{thm-01} and \ref{thm-02}  will be given in Sec. \ref{csw-sec3}.
%Theorems \ref{thm-01} and %\ref{thm-02} will be proved in the last Sec. \ref{csw-sec4}
%\section{Conjugate type   properties  of harmonic  $(K,K_{0})$-quasiregular mappings}\label{csw-sec2}

%The proofs of Proposition \ref{Prop-1} and Theorems \ref{thm-0.3} and \ref{eq-thm-9},  will be provided in Section \ref{csw-sec3}.
%Theorems \ref{thm-01} and \ref{thm-02} will be proved in the last Sec. \ref{csw-sec4}.
The proofs of Proposition \ref{Prop-1} and Theorems \ref{thm-0.3} and \ref{eq-thm-9} are presented in Section \ref{csw-sec3}, while Theorems \ref{thm-01} and \ref{thm-02} are proved in the final section, Section \ref{csw-sec4}.

\section{Radial growth type theorems   of harmonic  $(K,K_{0})$-quasiregular mappings}\label{csw-sec3}

%\section{Conjugate functions type   growth  theorems of harmonic quasiregular mappings}\label{csw-sec3}

\subsection*{The proof of Proposition \ref{Prop-1}} Since $f$ is a solution to (\ref{Bel-q}), we see that, for $z\in\Omega$,
\be\label{eq-ghl-1} |f_{\overline{z}}(z)|\leq\|\kappa_{1}\|_{\infty}|f_{z}(z)|+\|\kappa_{2}\|_{\infty}.\ee
It follows from (\ref{eq-ghl-1}) that
\be\label{eq-ghl-2}
\Lambda_{f}(z)\leq \sigma_{1}\lambda_{f}(z)+\sigma_{2}
\leq\sigma_{1}\lambda_{f}(z)+
\sqrt{\sigma_{1}^{2}(\lambda_{f}(z))^{2}+\sigma_{2}^{2}},
\ee
where $\sigma_{1}=(1+\|\kappa_{1}\|_{\infty})/(1-\|\kappa_{1}\|_{\infty})$ and $\sigma_{2}=2\|\kappa_{2}\|_{\infty}/(1-\|\kappa_{1}\|_{\infty})$.
The remaining proof will be divided into two cases.

\noindent $\mathbf{Case~1.}$ $\Lambda_{f}(z)\leq\sigma_{1}\lambda_{f}(z)$ for all $z\in\Omega$.

In this case, we have

$$\left(\Lambda_{f}(z)\right)^{2}\leq\,\sigma_{1}J_{f}(z),$$
which implies that $f$ is a $K$-quasiregular mapping, where $K=\sigma_{1}$.

\noindent $\mathbf{Case~2.}$ There is a subset $E$ of $\Omega$ such that
\beqq
\Lambda_{f}(z)>\sigma_{1}\lambda_{f}(z),~z\in E.
\eeqq

In this case, by (\ref{eq-ghl-2}), we have

\beqq
\left(\Lambda_{f}(z)-\sigma_{1}\lambda_{f}(z)\right)^{2}
\leq\sigma_{1}^{2}(\lambda_{f}(z))^{2}+\sigma_{2}^{2},~z\in\Omega.
\eeqq
This yields that
\be\label{eq-ghl-3}
(\Lambda_{f}(z))^{2}\leq2\sigma_{1}J_{f}(z)+\sigma_{2}^{2}.
\ee
Note that, for $z\in\Omega\setminus\,E$,
\be\label{eq-ghl-4}(\Lambda_{f}(z))^{2}\leq\sigma_{1}J_{f}(z).\ee
Combining (\ref{eq-ghl-3}) and (\ref{eq-ghl-4}) gives the desired result.
The proof of this Proposition is complete.
\qed

%For $p\in(0,\infty )$, denote by $L^{p}(\mathbb{T})$ the
%set of all measurable functions $\psi$ of $\mathbb{T}$ into
%$\mathbb{C}$ with
%$$\|\psi\|_{L^{p}}=\left(\frac{1}{2\pi}\int_{0}^{2\pi}|\psi(\zeta)|^{p}d\theta\right)^{\frac{1}{p}}<\infty.$$

 The following result is well-known.

\begin{Lem}\label{Lemx}
Suppose that $a,~b\in[0,\infty)$ and $p\in(0,\infty)$. Then
$$(a+b)^{p}\leq2^{\max\{p-1,0\}}(a^{p}+b^{p}).$$
\end{Lem}

Green's theorem  states that if $\varphi\in
\mathcal{C}^{2}(\mathbb{D})$, then
\be\label{eq1.2x}
\frac{1}{2\pi}\int_{0}^{2\pi}\varphi(re^{i\theta})\,d\theta=\varphi(0)+
\frac{1}{2}\int_{\mathbb{D}_{r}}\Delta (\varphi(z))\log\frac{r}{|z|}\,dA(z)
\ee  for $r\in (0, 1)$, where $dA=dxdy/\pi$  (see  \cite[Formula (2)]{P-09}).

Let $p\in(0,\infty)$, the spaces of Dirichlet type spaces $\mathcal{D}_{p-1}^{p}$ consists of all functions $f\in\mathscr{A}$
satisfying

$$\int_{\mathbb{D}}(1-|z|^{2})^{p-1}|f'(z)|^{p}dA(z)<\infty.$$
%where $dA(z)=dxdy/\pi$.
It follows from Littlewood-Paley's theorem that
 \be\label{rf-01}\mathcal{D}_{p-1}^{p}\subset\,H^{p}~\mbox{for}~ p\in(0,2],\ee
and  $H^{p}\subset\mathcal{D}_{p-1}^{p}$ for $p\in[2,\infty)$ (see \cite{G-P,L-P}).

%In particular, if $f$ is holomorphic in $\mathbb{D}$, then, for $p\in(0,\infty)$, we have
 %$$\frac{1}{2\pi}\int_{0}^{2\pi}|f(e^{i\theta})|^{p}\,d\theta=|f(0)|^{p}+
%\frac{1}{2}\int_{\mathbb{D}}\Delta \left(|f(z)|^{p}\right)\log\frac{1}{|z|}\,dA(z).$$

%\cite[Theorem~9.9]{gt}

\subsection*{The proof of Theorem \ref{thm-0.3}}
%  Since $\mathbb{D}$ is a simply connected domain, we see that $f$ admits a decomposition
% $f=h + \overline{g}$, where $h$ and $g$ are analytic
% in $\mathbb{D}$ with $g(0)=0$.
  %Let $F_{1}=h+g$ and $F_{2}=h-g$. Then $u={\rm Re}(F_{1})$ and $v={\rm Im}(F_{2})$.

 We first prove $(\mathscr{C}_{1})$.  Since $\mathbb{D}$ is a simply connected domain, we see that $f$ admits a decomposition
 $f=h + \overline{g}$, where $h$ and $g$ are analytic
 in $\mathbb{D}$ with $g(0)=0$. Let $F_{1}=h+g$ and $F_{2}=h-g$. Then $u={\rm Re}(F_{1})$ and $v={\rm Im}(F_{2})$.
 Since $f$ is a harmonic $(K,K_{0})$-quasiregular mapping in $\mathbb{D}$, we see that
$$(\Lambda_{f}(z))^{2}\leq KJ_{f}(z)+K_{0} ~\mbox{ for $z\in\mathbb{D}$}.
$$
This gives that, for $z\in\mathbb{D}$,
\beqq
\Lambda_{f}(z)\leq\frac{K\lambda_{f}(z)+\sqrt{\big(K\lambda_{f}(z))\big)^{2}+4K_{0}}}{2}
\leq K\lambda_{f}(z)+\sqrt{K_{0}},
\eeqq
and consequently
\be\label{K-1}|g'(z)|\leq\mu_{1}|h'(z)|+\mu_{2},
\ee where $\mu_{1}=(K-1)/(K+1)$ and $\mu_{2}=\sqrt{K_{0}}/(1+K)$.
By (\ref{K-1}), we have

\beqq
|F_{2}'|\leq\Lambda_{f}\leq(1+\mu_{1})|h'|+\mu_{2}
\eeqq
and
\beqq
|F_{1}'|\geq\lambda_{f}\geq(1-\mu_{1})|h'|-\mu_{2},
\eeqq
which imply that
\be\label{eq-ghk-2}
|F_{2}'|\leq\frac{1+\mu_{1}}{1-\mu_{1}}|F_{1}'|+\frac{2\mu_{2}}{1-\mu_{1}}=K|F_{1}'|+\sqrt{K_{0}}.
\ee

 From the assumption, we see that
\be\label{eq-jg-1}M_{p}(r,F_{1}')=M_{p}(r,\nabla u)=O\left(\frac{1}{1-r}\right),\ee 
 which, together with Theorem \ref{Thm-HL-3}, gives that
\be\label{eq-jg-2}M_{p}(r,F_{1})=O\left(\left(\log\frac{1}{1-r}\right)^{\frac{1}{p}}\right).\ee 
 It follows from (\ref{eq-ghk-2})  and (\ref{eq-jg-1}) that
 $$M_{p}(r,F_{2}')=O\left(\frac{1}{1-r}\right).$$ 
 By using Theorem \ref{Thm-HL-3} again, we have
 \be\label{eq-jg-3}M_{p}(r,F_{2})=O\left(\left(\log\frac{1}{1-r}\right)^{\frac{1}{p}}\right).\ee 
 Combing (\ref{eq-jg-2}) and (\ref{eq-jg-3}) gives $(\mathscr{C}_{1})$. Now we prove the sharpness part.
 Since all univalent analytic functions defined in $\mathbb{D}$ belong to $H^{p}$ for $p\in(0,1/2)$,  we assume $p\in[1/2,1)$ in the following extreme function.
 For $p\in[1/2,1)$, let
\be\label{edf01}f(z)=\frac{2K}{K+1}{\rm Re}(h_{0}(z))+i\frac{2}{K+1}{\rm Im}(h_{0}(z))=u(z)+iv(z),~z\in\mathbb{D},\ee where $h_{0}(z)=1/(1-z)^{1/p}$.
 It is not difficult to know that
 \be\label{lw-1}\frac{2}{K+1}|h_{0}|\leq|f|\leq\frac{2K}{K+1}|h_{0}|.\ee Then
 $$M_{p}(r,h_{0}')=M_{p}(r,\nabla u)=O\left(\frac{1}{1-r}\right).$$
  It follows from (\ref{lw-1}) and \cite[Remark 2]{G-P} that 
$$
M_p(r,f) \asymp \left(\log \frac{1}{1-r}\right)^{1/p}. 
$$ 
Here the symbol ``$\asymp$'' means that $M_p(r,f)$ and $\left(\log \frac{1}{1-r}\right)^{1/p}$ have the same order of growth as $r \to 1^-$.
%\[
%C_{1}^{\ast}\left(\log\frac{1}{1-r}\right)^{1/p}
%\le M_{p}(r,f)
%\le C_{2}^{\ast}\left(\log\frac{1}{1-r}\right)^{1/p}
%\]
%as $r\to 1^{-}$.}
  
  %$$M_{p}(r,f)=O\left(\left(\log\frac{1}{1-r}\right)^{\frac{1}{p}}\right)=M_{p}(r,h_{0})$$

 %By using  similar reasoning as in the proof of $(\mathscr{B}_{1})$, we can obtain $(\mathscr{B}_{2})$.

 %Then $\mathscr{F}_{r}\in\mathscr{D}_{p-1}^{p}$ for $p\in(1,2]$.

Next, we begin to prove $(\mathscr{C}_{2})$. For $z\in\mathbb{D}$, let $\mathscr{F}_{r}(z)=F_{1}(rz)=U_{r}(z)+iV_{r}(z),$ and
 let $\mathscr{F}_{r}^{\ast}(z)=F_{2}(rz)=U_{r}^{\ast}(z)+iV_{r}^{\ast}(z),$ where $r\in[0,1)$.
 Then $\mathscr{F}_{r},~\mathscr{F}_{r}^{\ast}\in\mathcal{D}_{p-1}^{p}$ for $p\in(1,2]$.
 By the closed graph theorem and (\ref{rf-01}), we
see that there is a positive constant $C=C(p)$ such that
\be\label{eq-GP-1-0}\|\mathscr{F}_{r}\|_{p}^{p}\leq\,C\left(|F_{1}(0)|^{p}+\mathscr{W}_{1}\right)<\infty,\ee
where $$\mathscr{W}_{1}=\int_{\mathbb{D}}(1-|z|^{2})^{p-1}|\mathscr{F}'_{r}(z)|^{p}dA(z).$$
Since
\be\label{eq-jg-x1}M_{p}(r,F_{1}')=M_{p}(r,\nabla u)=O\left(\eta(r)\right),\ee 
there is a positive constant $C=C(p)$ such that
\beqq
\mathscr{W}_{1}&=&\int_{\mathbb{D}}(1-|z|^{2})^{p-1}|\nabla U_{r}(z)|^{p}dA(z)\\
&=&2\int_{0}^{1}\rho(1-\rho^{2})^{p-1}\left(\int_{0}^{2\pi}|\nabla U_{r}(\rho\,e^{i\theta})|^{p}\frac{d\theta}{2\pi}\right)d\rho\\
&\leq&2C\int_{0}^{1}(1-\rho^{2})^{p-1}\big(\eta(r\rho)\big)^{p}d\rho\\
&=&2C\int_{0}^{r}(1-\rho^{2})^{p-1}\big(\eta(r\rho)\big)^{p}d\rho+2C\int_{r}^{1}(1-\rho^{2})^{p-1}\big(\eta(r\rho)\big)^{p}d\rho,
\eeqq
which, together with $r\rho\leq\,r$, $r\rho\leq\rho$ and the monotonic property of $\eta$, implies that
\beq\label{dx-0-1}
\mathscr{W}_{1}&\leq&2C\int_{0}^{r}(1-\rho^{2})^{p-1}\big(\eta(\rho)\big)^{p}d\rho+2^{p}C\big(\eta(r)\big)^{p}\int_{r}^{1}(1-\rho)^{p-1}d\rho\\ \nonumber
&=&2C\int_{0}^{r}(1-\rho^{2})^{p-1}\big(\eta(\rho)\big)^{p}d\rho+\frac{2^{p}C}{p}\big(\eta(r)\big)^{p}(1-r)^{p}.
\eeq

%for every $f\in\mathscr{D}_{p-1}^{p}$. %(see \cite[Theorem  B]{G-P}).
Since $$\|\mathscr{F}_{r}\|_{p}^{p}=\frac{1}{2\pi}\int_{0}^{2\pi}|F_{1}(re^{i\theta})|^{p}d\theta=M_{p}^{p}(r,F_{1}),$$
by (\ref{eq-GP-1-0}) and (\ref{dx-0-1}), we see that there is a positive constant $C=C(p)$ such that
\beq\label{dx-03}
M_{p}^{p}(r,u)&\leq&\,M_{p}^{p}(r,F_{1})\leq\,C\left(|F_{1}(0)|^{p}+\mathscr{W}_{1}\right)\\ \nonumber
&\leq&C\bigg[|F_{1}(0)|^{p}+2C\int_{0}^{r}(1-\rho^{2})^{p-1}\big(\eta(\rho)\big)^{p}d\rho\\  \nonumber
&&+\frac{2^{p}C}{p}\big(\eta(r)\big)^{p}(1-r)^{p}\bigg].
\eeq

%\beqq
%M_{p}^{p}(r,F_{1})&\leq&C\left(|F_{1}(0)|^{p}+\mathscr{W}_{1}\right)\\
%&=&C\left(|F_{1}(0)|^{p}+\int_{\mathbb{D}}(1-|z|^{2})^{p-1}|\nabla U_{r}(z)|^{p}dA(z)\right)\\
%&\leq&C\left(|F_{1}(0)|^{p}+\int_{0}^{1}(1-\rho^{2})^{p-1}\left(\int_{0}^{2\pi}|\nabla U_{r}(\rho\,e^{i\theta})|^{p}\frac{d\theta}{2\pi}\right)d\rho\right),
%\eeqq where $$\mathscr{W}_{1}=\int_{\mathbb{D}}(1-|z|^{2})^{p-1}|\mathscr{F}'_{r}(z)|^{p}dA(z).$$

%This together with (\ref{gfk-01}) gives
%\beqq
%M_{p}^{p}(r,F_{1})&\leq&C\left[|F_{1}(0)|^{p}+\int_{0}^{1}(1-\rho^{2})^{p-1}\big(\eta(r\rho)\big)^{p}d\rho\right]\\
%&=&C\left[|F_{1}(0)|^{p}+\int_{0}^{r}(1-\rho^{2})^{p-1}\big(\eta(r\rho)\big)^{p}d\rho+\int_{r}^{1}(1-\rho^{2})^{p-1}\big(\eta(r\rho)\big)^{p}d\rho\right]\\
%&\leq&C\left[|F_{1}(0)|^{p}+\int_{0}^{r}(1-\rho^{2})^{p-1}\big(\eta(\rho)\big)^{p}d\rho+2^{p-1}\big(\eta(r)\big)^{p}\int_{r}^{1}(1-\rho)^{p-1}d\rho\right]\\
%&=&C\left[|F_{1}(0)|^{p}+\int_{0}^{r}(1-\rho^{2})^{p-1}\big(\eta(\rho)\big)^{p}d\rho+\frac{2^{p-1}}{p}\big(\eta(r)\big)^{p}(1-r)^{p}\right].
%\eeqq

 On the other hand,
 using again the closed graph theorem and (\ref{rf-01}), we
see that there is a positive constant $C=C(p)$ such that
\be\label{eq-GP-1}\|\mathscr{F}_{r}^{\ast}\|_{p}^{p}\leq\,C\left(|F_{2}(0)|^{p}+\int_{\mathbb{D}}(1-|z|^{2})^{p-1}|(\mathscr{F}_{r}^{\ast}(z))'|^{p}dA(z)\right)<\infty.\ee
%Since $$\|\mathscr{F}_{r}^{\ast}\|_{p}^{p}=\frac{1}{2\pi}\int_{0}^{2\pi}|F_{2}(re^{i\theta})|^{p}d\theta=M_{p}^{p}(r,F_{2}),$$
From (\ref{eq-ghk-2}), (\ref{eq-GP-1}) and Lemma \ref{Lemx}, we see that there is a positive constant $C=C(p)$ such that

%which, together with (\ref{kp-0-1}) and Lemma \ref{Lemx}, yields that
\beq\label{eq-ad-1}
M_{p}^{p}(r,F_{2})&\leq&C\left(|F_{2}(0)|^{p}+\int_{\mathbb{D}}(1-|z|^{2})^{p-1}|\big(\mathscr{F}_{r}^{\ast}(z)\big)'|^{p}dA(z)\right)\\ \nonumber
&\leq&C\left(|F_{2}(0)|^{p}+\int_{\mathbb{D}}(1-|z|^{2})^{p-1}\big(K|\mathscr{F}'_{r}(z)|+\sqrt{K_{0}}\big)^{p}dA(z)\right)\\ \nonumber
&\leq&C\left(|F_{2}(0)|^{p}+2^{p-1}K^{p}\mathscr{W}_{1}+\mathscr{W}_{2}\right),
\eeq
where %$$\mathscr{W}_{1}=\int_{\mathbb{D}}(1-|z|^{2})^{p-1}|\mathscr{F}'_{r}(z)|^{p}dA(z)$$ and
 $$\mathscr{W}_{2}=2^{p-1}(K_{0})^{\frac{p}{2}}\int_{\mathbb{D}}(1-|z|^{2})^{p-1}dA(z)<\infty.$$
%Now, we come to estimate $\mathscr{W}_{1}$. By (\ref{gfk-01}), we have
%\beqq
%\mathscr{W}_{1}&=&\int_{\mathbb{D}}(1-|z|^{2})^{p-1}|\nabla U_{r}(z)|^{p}dA(z)\\
%&=&\int_{0}^{1}(1-\rho^{2})^{p-1}\left(\int_{0}^{2\pi}|\nabla U_{r}(\rho\,e^{i\theta})|^{p}\frac{d\theta}{2\pi}\right)d\rho\\
%&\leq&\int_{0}^{1}(1-\rho^{2})^{p-1}\big(\eta(r\rho)\big)^{p}d\rho\\
%&=&\int_{0}^{r}(1-\rho^{2})^{p-1}\big(\eta(r\rho)\big)^{p}d\rho+\int_{r}^{1}(1-\rho^{2})^{p-1}\big(\eta(r\rho)\big)^{p}d\rho,
%\eeqq
%which, together with $r\rho\leq\,r$, $r\rho\leq\rho$ and the monotonic property of $\eta$, implies that
%\beq\label{dx-01}
%\mathscr{W}_{1}&\leq&\int_{0}^{r}(1-\rho^{2})^{p-1}\big(\eta(\rho)\big)^{p}d\rho+2^{p-1}\big(\eta(r)\big)^{p}\int_{r}^{1}(1-\rho)^{p-1}d\rho\\ \nonumber
%&=&\int_{0}^{r}(1-\rho^{2})^{p-1}\big(\eta(\rho)\big)^{p}d\rho+\frac{2^{p-1}}{p}\big(\eta(r)\big)^{p}(1-r)^{p}.
%\eeq
Combing (\ref{dx-0-1}) and (\ref{eq-ad-1})  gives
\beq\label{jjl-1}
M_{p}(r,v)&\leq&\,M_{p}(r,F_{2})\\ \nonumber
&=&O\left(\left(2\int_{0}^{r}(1-\rho^{2})^{p-1}\big(\eta(\rho)\big)^{p}d\rho
+\frac{2^{p}}{p}\big(\eta(r)\big)^{p}(1-r)^{p}\right)^{\frac{1}{p}}\right)\\ \nonumber
&=&O\left(\left(\int_{0}^{r}(1-\rho^{2})^{p-1}\big(\eta(\rho)\big)^{p}d\rho
+\frac{2^{p-1}}{p}\big(\eta(r)\big)^{p}(1-r)^{p}\right)^{\frac{1}{p}}\right).
\eeq
Hence (\ref{Sh-11}) follows from (\ref{dx-03}) and (\ref{jjl-1}).

Now we  prove  (\ref{eq-101})  as follows:
Let $\rho\in(0,1)$. Then,
for $z=re^{it}\in\mathbb{D}_{\rho}$, we have
\beqq
f(z)=\frac{1}{2\pi}\int_{0}^{2\pi}\left(\frac{\rho e^{i(\theta+t)}}{\rho e^{i(\theta+t)}-z}+\frac{\overline{z}}{\rho e^{-i(\theta+t)}-\overline{z}}\right)f(\rho e^{i(\theta+t)})d\theta=
h(z)+\overline{g(z)}.
\eeqq
Then $$h'(z)=f_{z}(z)=\frac{1}{2\pi}\int_{0}^{2\pi}\frac{\rho e^{i(\theta+t)}}{(\rho e^{i(\theta+t)}-z)^{2}}f(\rho e^{i(\theta+t)})d\theta$$
and
$$\overline{g'(z)}=f_{\overline{z}}(z)=\frac{1}{2\pi}\int_{0}^{2\pi}\frac{\rho e^{-i(\theta+t)}}{(\rho e^{-i(\theta+t)}-\overline{z})^{2}}f(\rho e^{i(\theta+t)})d\theta,$$
which imply that

\beq\label{er-107}
|F'_{1}(z)|=|h'(z)+g'(z)|&=&\left|\frac{1}{2\pi}\int_{0}^{2\pi}\frac{\rho e^{i(\theta+t)}}{(\rho e^{i(\theta+t)}-z)^{2}}2{\rm Re}(f(\rho e^{i(\theta+t)}))d\theta\right|\\
\nonumber
&\leq&\frac{\rho}{\pi}\int_{0}^{2\pi}\frac{|f(\rho e^{i(\theta+t)})|}{|\rho e^{i\theta}-r|^{2}}d\theta.
\eeq
Now we set $\rho=(1+r)/2$.
It follows from (\ref{er-107}) and the Minkowski inequality that
\beqq
M_{p}(r,\nabla u)&=&M_{p}(r,F_{1}')\leq2\int_{0}^{2\pi}\left(\int_{0}^{2\pi}\frac{|f(\rho e^{i(\theta+t)})|^{p}}{|\rho e^{i\theta}-r|^{2p}}\frac{d\theta}{2\pi}\right)^{\frac{1}{p}}\frac{dt}{2\pi}\\
&=&\frac{M_{p}(\rho,f)}{\rho^{2}-r^{2}}\\
&=&O\left(\frac{\Psi(r)}{1-r}\right).
\eeqq

%The sharpness part of (\ref{103-eq})  follows from the following extreme function
%$$f(z)=\frac{2K}{K+1}{\rm Re}(h_{0}^{\ast}(z))+i\frac{2}{K+1}{\rm Im}(h_{0}^{\ast}(z)),~z\in\mathbb{D},$$ where
%$h_{0}^{\ast}(z)=1/(1-z)^{1/p}$ for $p\in(1,2)$.

At last, we prove $(\mathscr{C}_{3})$. By (\ref{eq1.2x}), we have
\beq\label{eq-jh-01}
M_{p}^{p}(r,F_{1})&=&|F_{1}(0)|^{p}+\frac{p^{2}}{2}\int_{\mathbb{D}_{r}}|F_{1}(z)|^{p-2}|F'_{1}(z)|^{2}\log\frac{r}{|z|}dA(z)\\ \nonumber
&=&|F_{1}(0)|^{p}+p^{2}\int_{0}^{r}\rho\log\frac{r}{\rho}\mathscr{J}(\rho)d\rho,
\eeq where
$$\mathscr{J}(\rho)=\frac{1}{2\pi}\int_{0}^{2\pi}|F_{1}(\rho e^{i\theta})|^{p-2}|F'_{1}(\rho e^{i\theta})|^{2}d\theta.$$
Applying H\"older's inequality, we have
\beqq
\mathscr{J}(\rho)\leq\,M_{p}^{p-2}(\rho,F_{1})M_{p}^{2}(\rho,F_{1}'),
\eeqq
which, together with (\ref{eq-jh-01}) and the fact ``$M_{p}(\rho, F_{1})$ being an increasing function of $\rho\in(0,r)$", implies that
\beq\label{e-fj-1}
M_{p}^{p}(r,F_{1})&\leq&|F_{1}(0)|^{p}+p^{2}\int_{0}^{r}
\left(\rho\log\frac{r}{\rho}\right)M_{p}^{p-2}(\rho,F_{1})M_{p}^{2}(\rho,F_{1}')d\rho\\ \nonumber
&\leq&|F_{1}(0)|^{p}+p^{2}M_{p}^{p-2}(r,F_{1})\int_{0}^{r}
\left(\rho\log\frac{r}{\rho}\right)M_{p}^{2}(\rho,F_{1}')d\rho.
\eeq
 Then, by (\ref{eq-jg-x1}) and (\ref{e-fj-1}), there is a positive constant $C$, independent of $F_{1}$, such that
 \beq\label{e-fj-2}
 M_{p}^{2}(r,F_{1})&\leq&|F_{1}(0)|^{2}+p^{2}\int_{0}^{r}
\left(\rho\log\frac{r}{\rho}\right)M_{p}^{2}(\rho,F_{1}')d\rho\\ \nonumber
&\leq&|F_{1}(0)|^{2}+p^{2}C\int_{0}^{r}\left(\rho\log\frac{r}{\rho}\right)(\eta(\rho))^{2}d\rho.
 \eeq
Since $$\rho\log\frac{r}{\rho}\leq\,r-\rho\leq1-\rho,$$ by (\ref{e-fj-2}), we see that
\be\label{e-fj-3}
M_{p}(r,u)\leq\,M_{p}(r,F_{1})\leq\left[|F_{1}(0)|^{2}+p^{2}C\int_{0}^{r}\left(1-\rho\right)(\eta(\rho))^{2}d\rho\right]^{\frac{1}{2}}
\ee

On the other hand, by (\ref{eq-ghk-2}) and (\ref{eq-jg-x1}), we have
 $$M_{p}(r,F_{2}')=O\left(\eta(r)\right).$$
 By using  similar reasoning as in the proof of (\ref{e-fj-3}), we obtain

\be\label{e-fj-4}
M_{p}(r,v)\leq\,M_{p}(r,F_{2})=O\left(\left(\int_{0}^{r}\left(1-\rho\right)(\eta(\rho))^{2}d\rho\right)^{\frac{1}{2}}\right).
\ee 
Therefore,  it follows from (\ref{e-fj-3}), (\ref{e-fj-4}) and  Minkowski's inequality that
\beqq
 M_{p}(r,f)&=&
 \left(\frac{1}{2\pi}\int_{0}^{2\pi}\left(|u(re^{i\theta})|^{2}+|v(re^{i\theta})|^{2}\right)^{\frac{p}{2}}d\theta\right)^{\frac{1}{p}}
 \leq\left(M_{p}^{2}(r,u)+M_{p}^{2}(r,v)\right)^{\frac{1}{2}}\\
 &\leq&M_{p}(r,u)+M_{p}(r,v)\\
 &=&O\left(\left(\int_{0}^{r}\left(1-\rho\right)(\eta(\rho))^{2}d\rho\right)^{\frac{1}{2}}\right).
\eeqq

By using similar reasoning as in the proof of (\ref{eq-101}), we can obtain
 (\ref{eq-102}).
The proof of this theorem is complete.
\qed

\subsection*{The proof of Corollary \ref{cor-2.2}} 
We first prove the sharpness part of (\ref{103-eq}).
Let $$f(z)=\frac{2K}{K+1}{\rm Re}(h_{0}^{\ast}(z))+i\frac{2}{K+1}{\rm Im}(h_{0}^{\ast}(z)),~z\in\mathbb{D},$$ where
$h_{0}^{\ast}(z)=1/(1-z)^{1/p}$ for $p\in(1,2)$.
Since
$$\frac{2}{K+1}|h_{0}^{\ast}|\leq|f|\leq\frac{2K}{K+1}|h_{0}^{\ast}|,$$
it follows from \cite[Remark 2]{G-P} that $$
M_p(r,f) \asymp \left(\log \frac{1}{1-r}\right)^{1/p}.
$$

In order to prove the sharpness part of (\ref{104-eq}), let $$f(z)=\frac{2K}{K+1}{\rm Re}(h^{\ast}(z))+i\frac{2}{K+1}{\rm Im}(h^{\ast}(z))=u(z)+iv(z),~z\in\mathbb{D},$$
where $h^{\ast}(z)=\sum_{n=0}^{\infty}z^{2^{n}}$. It follows from \cite[Lemma~2.1]{A-C-P} that
$$\frac{K+1}{2K}M_{p}(r,\nabla u)=M_{p}(r,(h^{\ast})')=O\left(\frac{1}{1-r}\right).
$$ By \cite[Theorem~A]{G-P}, we have
\be\label{qd-01}
M_{p}(r,h^{\ast})=O\left(\left(\log\frac{1}{1-r}\right)^{\frac{1}{2}}\right).
\ee 
Since $$\frac{2}{K+1}|h^{\ast}|\leq|f|\leq\frac{2K}{K+1}|h^{\ast}|,$$ by
(\ref{qd-01}), we see that
$$M_{p}(r,f)\asymp O\left(\left(\log\frac{1}{1-r}\right)^{\frac{1}{2}}\right).$$
The proof of this corollary is complete.
\qed

%\begin{Thm}{\rm  (\cite[Theorem  2.11]{CH-2023-SCM})}\label{Thm-g5.0}
%Let $\omega\in E$, and let $p\in[1,\infty]$ be a constant. %Then the following statements are equivalent.

%\begin{enumerate}
%\item[{\rm (a)}] If $\psi\in\mathscr{L}_{\omega,p}(\mathbb{T})$ and
%$\psi$ is continuous on $\mathbb{T}$, then $P[\psi]\in\mathscr{L}_{\omega,p}(\overline{\mathbb{D}})$;
%\item[{\rm (b)}] There is a positive constant $C$ such that for all $\delta\in[0,\pi]$,
%$$\delta\int_{\delta}^{\pi}\frac{\omega(t)}{t^{2}}dt\leq\,C\omega(\delta).$$
%\end{enumerate}
%Then $(b)\Rightarrow(a)$ for $p\in[1,\infty]$, and $(a)\Leftrightarrow(b)$ for $p=\infty$.
%\end{Thm}

\subsection*{The proof of Theorem \ref{eq-thm-9}} 
We  divide the proof of this theorem into two cases.

\noindent  $\mathbf{Case~1}$. Let $p\in[1,\infty)$.

%Without loss of generality, we assume $p\in[1,\infty)$.
 %Since $\mathscr{D}_{2}$ can be deduced by $\mathscr{D}_{1}$,
%we only need to prove $\mathscr{D}_{1}$.
%We divide the proof of this theorem into two cases.
We first prove the necessity. It follows from the concavity of $\varphi\in E$ that $\varphi(t)/t$ is decreasing for $t\in(0,\infty)$.
Since $u$ is continuous in $\overline{\mathbb{D}}$ and $u\in\mathscr{L}_{\varphi,p}(\mathbb{T})$,
by \cite[Theorem  2.11]{CH-2023-SCM}, we see that $u\in\mathscr{L}_{\varphi,p}(\overline{\mathbb{D}})$.
 For any fixed $z\in\mathbb{D}$, let $r=d_{\mathbb{D}}(z)/2$, where $d_{\mathbb{D}}(z)=1-|z|$. Then
there is a positive constant $C$ such that
\be\label{eq-0.01k}
\left(\int_{0}^{2\pi}|u(e^{i\theta}z)-u(e^{i\theta}\xi)|^{p}d\theta\right)^{\frac{1}{p}}\leq C\varphi(|z-\xi|)\leq C\varphi(r)
\ee
for $\xi\in\overline{\mathbb{D}(z,r)}$.
For $\theta\in[0,2\pi]$ and $e^{i\theta}\xi\in\mathbb{D}(e^{i\theta}z,r)$, we have
$$u(e^{i\theta}\xi)=\frac{1}{2\pi}\int_{0}^{2\pi}P_{r}(\xi,e^{i\eta})u(ze^{i\theta}+e^{i\theta}re^{i\eta})d\eta,$$
where $$P_{r}(\xi,e^{i\eta})=\frac{r^{2}-|\xi-z|^{2}}{|re^{i\eta}-(\xi-z)|^{2}}.$$
Elementary calculations lead to
\beqq
\frac{\partial}{\partial\xi}P_{r}(\xi,e^{i\eta})&=&\frac{\overline{z}-\overline{\xi}}{|re^{i\eta}-(\xi-z)|^{2}}\\
&&+
\frac{(r^{2}-|\xi-z|^{2})(re^{-i\eta}-(\overline{\xi}-\overline{z}))}{|re^{i\eta}-(\xi-z)|^{4}}
\eeqq
and
\beqq
\frac{\partial}{\partial\overline{\xi}}P_{r}(\xi,e^{i\eta})&=&\frac{z-\xi}{|re^{i\eta}-(\xi-z)|^{2}}\\
&&+
\frac{(r^{2}-|\xi-z|^{2})(re^{i\eta}-(\xi-z))}{|re^{i\eta}-(\xi-z)|^{4}}.
\eeqq
Then, for $e^{i\theta}\xi\in\mathbb{D}(e^{i\theta}z,r/2)$, we have
\be\label{eq-0.1k}
\left|\frac{\partial}{\partial\xi}P_{r}(\xi,e^{i\eta})\right|\leq\frac{10}{r}
\ee
and
\be\label{eq-0.2k}
\left|\frac{\partial}{\partial\overline{\xi}}P_{r}(\xi,e^{i\eta})\right|\leq\frac{10}{r}.
\ee
It follows from (\ref{eq-0.1k}),  (\ref{eq-0.2k}) and H\"{o}lder's inequality that
\beq\label{eq-0.3k}
\big(\Lambda_{u}(e^{i\theta}\xi)\big)^{p}&\leq&\int_{0}^{2\pi}\left(\Lambda_{P_{r}}(\xi, e^{i\eta})\right)^p
\left|u(ze^{i\theta}+e^{i\theta}re^{{\rm i}\eta})-u(ze^{i\theta})\right|^{p}\frac{d\eta}{2\pi}\\ \nonumber
&\leq&\frac{(20)^{p}}{r^{p}}\int_{0}^{2\pi}
\left|u(ze^{i\theta}+e^{i\theta}re^{{\rm i}\eta})-u(ze^{i\theta})\right|^{p}\frac{d\eta}{2\pi}.
\eeq %where $\mathscr{P}_{f}(\theta,\eta)=f(ze^{{\rm i}\theta}+e^{{\rm i}\theta}re^{{\rm i}\eta})-f(ze^{{\rm i}\theta})$.
By taking $\xi=z$ in (\ref{eq-0.3k}) and integrating both sides of the inequality (\ref{eq-0.3k}) with respect to $\theta$ from $0$ to $2\pi$, we obtain  from (\ref{eq-0.01k}) that
\beq\label{ed-r1}
\Gamma_{u}(z)&:=&\left(\int_{0}^{2\pi}\left(\Lambda_{u}(ze^{i\theta})\right)^{p}d\theta\right)^{\frac{1}{p}}\\  \nonumber
&\leq&
\frac{20}{r}\left(\int_{0}^{2\pi}\left(\int_{0}^{2\pi}
\left|u(ze^{i\theta}+e^{i\theta}re^{i\eta})-u(ze^{i\theta})\right|^{p}\frac{d\eta}{2\pi}\right) d\theta\right)^{\frac{1}{p}}
\\ \nonumber
&=&
\frac{20}{r}\left(\int_{0}^{2\pi}\left(\int_{0}^{2\pi}
\left|u(ze^{i\theta}+e^{i\theta}re^{{\rm i}\eta})-u(ze^{i\theta})\right|^{p}d\theta\right) \frac{d\eta}{2\pi}\right)^{\frac{1}{p}}
\\ \nonumber
&\leq&
40C\frac{\varphi\left(\frac{d_{\mathbb{D}}(z)}{2}\right)}{d_{\mathbb{D}}(z)}\\ \nonumber
 &\leq& 40C\frac{\varphi\left(d_{\mathbb{D}}(z)\right)}{d_{\mathbb{D}}(z)}.
 \eeq %where $$\Gamma_{u}(z)=\left(\int_{0}^{2\pi}\left(\Lambda_{u}(ze^{i\theta})\right)^{p}d\theta\right)^{\frac{1}{p}}.$$
% Since $\omega(t)/t$ is nonincreasing for $t\in(0,\infty)$, we see that
 % $$\frac{\omega\left(d_{\mathbb{D}}(z)\right)}{d_{\mathbb{D}}(z)}\leq
  %\frac{\omega\left(\frac{d_{\mathbb{D}}(z)}{2}\right)}{\frac{d_{\mathbb{D}}(z)}{2}}$$
 Similar to the proof of Theorem \ref{thm-0.3}, we let $f=h + \overline{g}$, where $h$ and $g$ are analytic
 in $\mathbb{D}$ with $g(0)=0$. Let $F_{1}=h+g$ and $F_{2}=h-g$. Then $u={\rm Re}(F_{1})$ and $v={\rm Im}(F_{2})$.
Since $|\nabla u|=|F_{1}'|$ and $|\nabla v|=|F_{2}'|$, by (\ref{eq-ghk-2}) and Lemma \ref{Lemx},
we see that
\beq\label{ed-r2}
\Lambda_{v}^{p}&=&|\nabla v|^{p}\leq\left(K|\nabla u|+\sqrt{K_{0}}\right)^{p}\leq2^{p-1}\left(K^{p}|\nabla u|^{p}+(K_{0})^{\frac{p}{2}}\right)\\ \nonumber
&=&2^{p-1}\left(K^{p}\Lambda_{u}^{p}+(K_{0})^{\frac{p}{2}}\right).
\eeq
By Lemma \ref{Lemx}, we have
\beqq
\int_{0}^{2\pi}\left(\Lambda_{f}(ze^{i\theta})\right)^{p}d\theta&\leq&\int_{0}^{2\pi}\left(\Lambda_{u}(ze^{i\theta})+\Lambda_{v}(ze^{i\theta})\right)^{p}d\theta\\
&\leq&2^{p-1}(\Gamma_{u}(z))^{p}+2^{p-1}\int_{0}^{2\pi}\left(\Lambda_{v}(ze^{i\theta})\right)^{p}d\theta,
\eeqq
which, together with (\ref{ed-r1}) and (\ref{ed-r2}), implies that
\beq\label{ed-r3}
\left(\int_{0}^{2\pi}\left(\Lambda_{f}(ze^{i\theta})\right)^{p}d\theta\right)^{\frac{1}{p}}&\leq&
\left(\left(2^{p-1}+2^{2p-2}K^{p}\right)(\Gamma_{u}(z))^{p}+\pi2^{2p-1}(K_{0})^{\frac{p}{2}}\right)^{\frac{1}{p}}\\ \nonumber
&\leq&\left(2^{p-1}+2^{2p-2}K^{p}\right)^{\frac{1}{p}}\Gamma_{u}(z)+\pi^{\frac{1}{p}}2^{\frac{2p-1}{p}}\sqrt{K_{0}}\\ \nonumber
&\leq&40C\left(2^{p-1}+2^{2p-2}K^{p}\right)^{\frac{1}{p}}\frac{\varphi\big(d_{\mathbb{D}}(z)\big)}{d_{\mathbb{D}}(z)}\\ \nonumber
&&+\pi^{\frac{1}{p}}2^{\frac{2p-1}{p}}\sqrt{K_{0}},
\eeq
where $C$ is a positive constant.
%Since $\varphi(t)/t$ is nonincreasing for $t\in(0,\infty)$,
%we see that
%$$\frac{\varphi\big(d_{\mathbb{D}}(z)\big)}{d_{\mathbb{D}}(z)}\geq\varphi(1),$$ which, together with
% (\ref{ed-r3}), implies that there is a positive constant $C$ such that
%\beqq
%\left(\int_{0}^{2\pi}\left(\Lambda_{f}(ze^{i\theta})\right)^{p}d\theta\right)^{\frac{1}{p}}\leq\,C\frac{\varphi\big(d_{\mathbb{D}}(z)\big)}{d_{\mathbb{D}}(z)}.
%\eeqq

Next, we prove the sufficiency. 

For $\xi\in\overline{\mathbb{D}}$, let $f_{R}(\xi)=f(R\xi)=U_{R}(\xi)+iV_{R}(\xi)=u(R\xi)+iv(R\xi)$, where $R\in(0,1)$.
 Given $\zeta_{1},~\zeta_{2}\in\mathbb{T}$, let $\varrho=|\zeta_{1}-\zeta_{2}|/2$, and let $L_{0}$ be the shorter subarc of $\mathbb{T}$
connecting $\zeta_{1}$ and $\zeta_{2}$. Set $\gamma_{1}=\{r\zeta_{1}:~r\in[1-\varrho, 1]\}$, $\gamma_{2}=\{(1-\varrho)\zeta:~\zeta\in L_{0}\}$,
$\gamma_{3}=\{r\zeta_{2}:~r\in[1-\varrho,1]\}$, and $\gamma=\gamma_{1}\cup\gamma_{2}\cup\gamma_{3}$. When endowed with the appropriate orientation, $\gamma$
becomes a path connecting  $\zeta_{1}$ to $\zeta_{2}$. %Then we have
%\beqq
%|u(\zeta_{1})-u(\zeta_{2})|\leq\int_{\gamma}|\nabla\,u(z)||dz|\leq\int_{\gamma}\Lambda_{f}(z)|dz|.
%\eeqq
By Minkowski's inequality, we have

\beq\label{ed-r6}
\mathcal{L}_{p}[U_{R}](\zeta_{1},\zeta_{2})&\leq&
\left(\int_{0}^{2\pi}\left(\int_{\gamma}|\nabla U_{R}(ze^{{\rm i}\theta})|ds(z)\right)^{p}d\theta\right)^{\frac{1}{p}}\\ \nonumber
&\leq&
\int_{\gamma}\left(\int_{0}^{2\pi}|\nabla U_{R}(ze^{{\rm i}\theta})|^{p}d\theta\right)^{\frac{1}{p}}ds(z),
\eeq where $ds$ stands for the arc length measure on
$\gamma$.
Since $|\nabla U_{R}|\leq\Lambda_{f_{R}}$ by (\ref{gv-1}) and (\ref{ed-r6}), we see that there are positive constants $C(p)$ and $C(K,p)$ such that

\beq\label{ed-r7}
\mathcal{L}_{p}[U_{R}](\zeta_{1},\zeta_{2})&\leq&\int_{\gamma}\left(\int_{0}^{2\pi}\left(\Lambda_{f_{R}}(ze^{{\rm i}\theta})\right)^{p}d\theta\right)^{\frac{1}{p}}ds(z)\\ \nonumber
&=&R\int_{\gamma}\left(\int_{0}^{2\pi}\left(\Lambda_{f}(Rze^{{\rm i}\theta})\right)^{p}d\theta\right)^{\frac{1}{p}}ds(z)\\ \nonumber
&\leq&C(K,p)\int_{\gamma}\frac{\varphi\big(d_{\mathbb{D}}(Rz)\big)}{d_{\mathbb{D}}(Rz)}ds(z)+C(p)\sqrt{K_{0}}\int_{\gamma}ds(z)\\ \nonumber
&=&C(K,p)\sum_{j=1}^{3}\int_{\gamma_{j}}\frac{\varphi\big(d_{\mathbb{D}}(Rz)\big)}{d_{\mathbb{D}}(Rz)}ds(z)+C(p)\sqrt{K_{0}}\sum_{j=1}^{3}\int_{\gamma_{j}}ds(z).
\eeq
Since $\varphi(t)/t$ is nonincreasing for $t\in(0,\infty)$, it follows from (\ref{eq2x}) that there is a positive constant $C$ such that
\be\label{ed-r8}
\int_{\gamma_{1}}\frac{\varphi\big(d_{\mathbb{D}}(Rz)\big)}{d_{\mathbb{D}}(Rz)}ds(z)\leq\int_{\gamma_{1}}\frac{\varphi\big(d_{\mathbb{D}}(z)\big)}{d_{\mathbb{D}}(z)}ds(z)
\leq\,C\varphi(R|\zeta_{1}-\zeta_{2}|)
\ee
and
\be\label{ed-r9}
\int_{\gamma_{3}}\frac{\varphi\big(d_{\mathbb{D}}(Rz)\big)}{d_{\mathbb{D}}(Rz)}ds(z)\leq\int_{\gamma_{3}}\frac{\varphi\big(d_{\mathbb{D}}(z)\big)}{d_{\mathbb{D}}(z)}ds(z) \leq\,C\varphi(|R\zeta_{1}-R\zeta_{2}|).
\ee
On the other hand,
\beq\label{ed-r10}
\int_{\gamma_{2}}\frac{\varphi\big(d_{\mathbb{D}}(Rz)\big)}{d_{\mathbb{D}}(Rz)}ds(z)&\leq&\int_{\gamma_{2}}\frac{\varphi\big(d_{\mathbb{D}}(z)\big)}{d_{\mathbb{D}}(z)}ds(z)\\ \nonumber
&=&\frac{\varphi\big(\varrho\big)}{\varrho}\int_{\gamma_{2}}ds(z)\\ \nonumber
&=&\frac{2\arcsin\varrho}{\varrho}\varphi\big(\varrho\big)\\ \nonumber
&\leq&\,C_{0}\varphi(|\zeta_{1}-\zeta_{2}|),
\eeq where $$C_{0}=\max_{\varrho\in[0,1]}\left\{\frac{2\arcsin\varrho}{\varrho}\right\}.$$

Since $\varphi(t)/t$ is nonincreasing for $t\in(0,\infty)$, we see that
$$\frac{\varphi(|\zeta_{1}-\zeta_{2}|)}{|\zeta_{1}-\zeta_{2}|}\geq\frac{\varphi(2)}{2}$$
which implies that
there is a positive $C$ such that
\be\label{ed-r1-9}\sum_{j=1}^{3}\int_{\gamma_{j}}ds(z)\leq C|\zeta_{1}-\zeta_{2}|\leq\frac{2C}{\varphi(2)}\varphi(|\zeta_{1}-\zeta_{2}|).\ee
Combining (\ref{ed-r7}), (\ref{ed-r8}), (\ref{ed-r9}), (\ref{ed-r10}) and (\ref{ed-r1-9}) gives $U_{R}\in\mathscr{L}_{\varphi,p}(\mathbb{T})$.
This, together with the continuity of $u$  in $\overline{\mathbb{D}}$, yields that there exists  a positive constant $C$
such that, for all $\xi_{1},~\xi_{2}\in\mathbb{T}$,
 $$\mathcal{L}_{p}[u](\xi_{1},\xi_{2})=\lim_{R\rightarrow1^{-}}\mathcal{L}_{p}[U_{R}](\xi_{1},\xi_{2})\leq C\varphi(|\xi_{1}-\xi_{2}|).$$
 
 \noindent  $\mathbf{Case~2}$. Let $p=\infty$.

The approach employed in proving Case 1 remains applicable to Case 2. It suffices to substitute the expressions involving
$p\in[1,\infty)$ with their counterparts for
$p=\infty$ to conclude the argument. Hence, the details are omitted here.
The proof of this theorem is complete. \qed

\section{Coefficient growth type theorems of univalent harmonic  $(K,K_{0})$-quasiregular mappings}\label{csw-sec4}
\subsection*{The proof of Theorem \ref{thm-01}}
We first prove (\ref{eq-th-c1}). Since  ${\mathcal
S}_{H}(K)\subset{\mathcal S}_{H}$, we see that

\be\label{eq-ry-1}\gamma^{\ast}\leq\gamma=\sup_{f\in{\mathcal
S}_{H}}\frac{|f_{zz}(0)|}{2}.\ee
Note that, for $h\in{\mathcal S}$, $$f=h+\frac{K-1}{K+1}\overline{h}\in{\mathcal
S}_{H}(K).$$ Then, by De Branges's result in \cite{B}, we have $$\gamma^{\ast}\geq2,$$ which, together with (\ref{eq-ry-1}), gives that
\be\label{eq-ry-2}2\leq\gamma^{\ast}\leq\gamma.\ee
%\begin{lem}
%Let $r\in(0,1)$, $p\in(0,\infty)$ and $f=h+\overline{g}\in{\mathcal S}_{H}(K)$, where $h$ and $g$ are analytic in $\mathbb{D}$. Then $$M_{p}(r,f)=O\left(\frac{1}{(1-r)^{\gamma-\frac{1}{p}}}\right)~\mbox{as}~r\rightarrow1^{-}.$$
%\end{lem}
%\bpf
For $f\in{\mathcal
S}_{H}(K)$ and any fixed $\zeta\in\mathbb{D}$, let
\be\label{hh-x1}F(z)=\frac{f\left(\frac{\zeta+z}{1+\overline{\zeta}z}\right)-f(\zeta)}{(1-|\zeta|^{2})h'(\zeta)}=H(z)+\overline{G(z)},~z\in\mathbb{D}.
\ee
Then $F\in{\mathcal S}_{H}(K)$ and
$$H(z)=z+A_{2}(\zeta)z^{2}+A_{3}(\zeta)z^{3}+\cdots,$$
where $$A_{2}(\zeta)=\frac{1}{2}\left[(1-|\zeta|^{2})\frac{h''(\zeta)}{h'(\zeta)}-2\overline{\zeta}\right].$$
Since $|A_{2}(\zeta)|\leq\gamma^{\ast}$, we see that
\be\label{hh-x2}{\rm  Re}\left\{\frac{\zeta h''(\zeta)}{h'(\zeta)}\right\}\leq\frac{2|\zeta|^{2}+2\gamma^{\ast}|\zeta|}{1-|\zeta|^{2}}.\ee
By (\ref{hh-x2}), we have
$$\frac{\partial}{\partial |\zeta|}\big(\log|h'(\zeta)|\big)\leq\frac{2|\zeta|+2\gamma^{\ast}}{1-|\zeta|^{2}},$$
which implies that
\be\label{hh-x3}
|h'(\zeta)|
\leq\frac{(1+|\zeta|)^{\gamma^{\ast}-1}}{(1-|\zeta|)^{\gamma^{\ast}+1}}.\ee

Note that $|g'|\leq\frac{K-1}{K+1}|h'|$, which, together with (\ref{hh-x3}), gives that
 \be\label{hh-x4}|g'(\zeta)|
\leq\frac{K-1}{K+1}\frac{(1+|\zeta|)^{\gamma^{\ast}-1}}{(1-|\zeta|)^{\gamma^{\ast}+1}}.\ee
It follows from (\ref{hh-x3}) and  (\ref{hh-x4}) that
\beq\label{hh-x5}
|f(z)|&=&\left|\int_{0}^{z}h'(\zeta)d\zeta+\overline{g'(\zeta)}d\overline{\zeta}\right|
\leq\frac{2K}{1+K}\int_{0}^{|z|}\frac{(1+t)^{\gamma^{\ast}-1}}{(1-t)^{\gamma^{\ast}+1}}dt\\ \nonumber
&=&\frac{K}{\gamma^{\ast}(1+K)}\left[\left(\frac{1+|z|}{1-|z|}\right)^{\gamma^{\ast}}-1\right].
\eeq

$\mathbf{Case~1}.$ Let $p\in(2,\infty)$.

By using   (\ref{eq1.2x}), we have
\beqq
M_{p}^{p}(r,f)&=&\frac{1}{2}\int_{\mathbb{D}_{r}}\Delta(|f(z)|^{p})\log\frac{r}{|z|}dA(z)\\
&=&\frac{1}{2}\int_{\mathbb{D}_{r}}\Big(p(p-2)|f(z)|^{p-4}\left|f(z)\overline{f_{z}(z)}+f_{\overline{z}}(z)\overline{f(z)}\right|^{2}\\
&&+
2p|f(z)|^{p-2}\left(|f_{z}(z)|^{2}+|f_{\overline{z}}(z)|^{2}\right)\Big)\log\frac{r}{|z|}dA(z),
\eeqq
which yields that

\beq\label{hh-x9}
r\frac{d}{dr}M_{p}^{p}(r,f)&\leq&\frac{p^{2}}{2}\int_{\mathbb{D}_{r}}|f(z)|^{p-2}\Lambda_{f}^{2}(z)dA(z)\\ \nonumber
&\leq&\frac{p^{2}K}{2}\int_{\mathbb{D}_{r}}|f(z)|^{p-2}J_{f}(z)dA(z)\\ \nonumber
&\leq&\frac{p^{2}K}{2}\int_{\mathbb{D}_{m(r)}}|w|^{p-2}dA(w),
\eeq where
\be\label{mr}m(r)=\max_{\theta\in[0,2\pi]}|f(re^{i\theta})|.\ee
By (\ref{eq-ry-2}), (\ref{hh-x5}) and (\ref{hh-x9}), we have
\beq\label{hh-x10}
M_{p}^{p}(r,f)&\leq& Kp\int_{0}^{r}\frac{(m(\rho))^{p}}{\rho}d\rho
\leq\frac{K^{2}p}{\gamma^{\ast}(1+K)}\int_{0}^{r}\frac{\left[\left(\frac{1+\rho}{1-\rho}\right)^{\gamma^{\ast}}-1\right]^{p}}{\rho}d\rho\\ \nonumber
&=&\frac{K^{2}p}{\gamma^{\ast}(1+K)}
\int_{0}^{r}\frac{[(1+\rho)^{\gamma^{\ast}}-(1-\rho)^{\gamma^{\ast}}]^{p}}{\rho}\frac{1}{(1-\rho)^{\gamma^{\ast} p}}d\rho.
\eeq
Since $$\lim_{\rho\rightarrow0^{+}}\frac{[(1+\rho)^{\gamma^{\ast}}-(1-\rho)^{\gamma^{\ast}}]^{p}}{\rho}=0,$$ we see that
$[(1+\rho)^{\gamma^{\ast}}-(1-\rho)^{\gamma^{\ast}}]^{p}/\rho$ is bounded for $\rho\in(0,1]$.
Let \be\label{tt-1}C_{1}=\max_{\rho\in(0,1]}\frac{[(1+\rho)^{\gamma^{\ast}}-(1-\rho)^{\gamma^{\ast}}]^{p}}{\rho}.\ee
Then, by (\ref{hh-x10}), we have

\beqq
M_{p}^{p}(r,f)&\leq&\frac{K^{2}pC_{1}}{\gamma^{\ast}(1+K)}\int_{0}^{r}\frac{1}{(1-\rho)^{\gamma^{\ast} p}}d\rho\\
&=&\frac{K^{2}pC_{1}}{\gamma^{\ast}(1+K)(p\gamma^{\ast} -1)}\left[\frac{1}{(1-r)^{p\gamma^{\ast} -1}}-1\right],
\eeqq
which implies that (\ref{eq-th-c1}) holds.

%\begin{lem}
%Let $f$ be a $(K,K_{0})$-quasiconformal and harmonic mapping of $\mathbb{D}$ into $\mathbb{C}$ with $f(0)=0$.
%Then, for $p\in(0,2]$, we have
%$$M_{p}^{p}(r,f)\leq2K\int_{0}^{r}\frac{(m(\rho))^{p}}{\rho}d\rho.$$
%\end{lem}
%\bpf

$\mathbf{Case~2}.$ Let $p\in(1/\gamma^{\ast},2]$.

For $n\in\{1,2,\ldots\}$, let $F_{n}(z)=\left(|f(z)|^{2}+\frac{1}{n}\right)^{\frac{1}{2}}$, $z\in\mathbb{D}$.
Elementary calculations yield that

\be\label{eq-1g}
\Delta(F_{n}^{p})=p(p-2)\left(|f|^{2}+\frac{1}{n}\right)^{\frac{p}{2}-2}\left|f\overline{f_{z}}+f_{\overline{z}}\overline{f}\right|^{2}+
2p\left(|f|^{2}+\frac{1}{n}\right)^{\frac{p}{2}-1}(|f_{z}|^{2}+|f_{\overline{z}}|^{2}).
\ee
It follows from    the Green Theorem that
$$M_{p}^{p}(r,F_{n})=(F_{n}(0))^{p}+\frac{1}{2}\int_{\mathbb{D}_{r}}\Delta(|F_{n}(z)|^{p})\log\frac{r}{|z|}dA(z),$$
which, together with (\ref{eq-1g}), implies that

\beq\label{eq-2g}
r\frac{d}{dr}M_{p}^{p}(r,F_{n})&=&\frac{1}{2}\int_{\mathbb{D}_{r}}\Delta(|F_{n}(z)|^{p})dA(z)\\ \nonumber
&\leq&p\int_{\mathbb{D}_{r}}\left(|f(z)|^{2}+\frac{1}{n}\right)^{\frac{p}{2}-1}\Lambda_{f}^{2}(z)dA(z)\\ \nonumber
&\leq&Kp\int_{\mathbb{D}_{r}}\left(|f(z)|^{2}+\frac{1}{n}\right)^{\frac{p}{2}-1}J_{f}(z)dA(z).
\eeq
  Then, by  (\ref{mr}) and (\ref{eq-2g}), we have

\beqq
r\frac{d}{dr}M_{p}^{p}(r,F_{n})&\leq&Kp\int_{\mathbb{D}_{m(r)}}\frac{d\sigma(w)}{\left(|w|^{2}+\frac{1}{n}\right)^{1-\frac{p}{2}}}\\ \nonumber
&=&\frac{Kp}{\pi}\int_{0}^{2\pi}\int_{0}^{m(r)}\frac{\rho}{\left(\rho^{2}+\frac{1}{n}\right)^{1-\frac{p}{2}}}d\rho d\theta\\ \nonumber
&=&2K\left(\left(m^{2}(r)+\frac{1}{n}\right)^{\frac{p}{2}}-\frac{1}{n^{\frac{p}{2}}}\right),
\eeqq
which, together with dominated convergence theorem, gives that
\be\label{eq-2} M_{p}^{p}(r,f)\leq2K\int_{0}^{r}\frac{(m(\rho))^{p}}{\rho}d\rho.\ee

Without loss of generality, we assume that $r>1/2$. Then, by (\ref{hh-x5}) and (\ref{eq-2}), we have

\beq\label{eq-2.0}
M_{p}^{p}(r,f)
&\leq&\frac{2K^{2}}{\gamma^{\ast}(1+K)}\int_{0}^{r}\frac{[(1+\rho)^{\gamma^{\ast}}-(1-\rho)^{\gamma^{\ast}}]^{p}}{\rho}\frac{1}{(1-\rho)^{\gamma^{\ast} p}}d\rho\\ \nonumber
&=&\frac{2K^{2}}{\gamma^{\ast}(1+K)}\int_{0}^{\frac{1}{2}}\frac{[(1+\rho)^{\gamma^{\ast}}-(1-\rho)^{\gamma^{\ast}}]^{p}}{\rho}\frac{1}{(1-\rho)^{\gamma^{\ast} p}}d\rho\\ \nonumber
&&+\frac{2K^{2}}{\gamma^{\ast}(1+K)}\int_{\frac{1}{2}}^{r}\frac{[(1+\rho)^{\gamma^{\ast}}-(1-\rho)^{\gamma^{\ast}}]^{p}}{\rho}\frac{1}{(1-\rho)^{\gamma^{\ast} p}}d\rho\\ \nonumber
&\leq&\frac{2^{1+p\gamma^{\ast}}K^{2}}{\gamma^{\ast}(1+K)}\varpi_{p}+
\frac{2K^{2}C_{2}}{\gamma^{\ast}(1+K)}\int_{\frac{1}{2}}^{r}\frac{1}{(1-\rho)^{\gamma^{\ast} p}}d\rho\\ \nonumber
&\leq&\frac{2^{1+p\gamma^{\ast}}K^{2}}{\gamma^{\ast}(1+K)}\varpi_{p}+
\frac{2K^{2}C_{2}}{\gamma^{\ast}(1+K)(p\gamma^{\ast}-1)}\left[(1-r)^{1-p\gamma^{\ast}}-\frac{1}{2^{1-p\gamma^{\ast}}}\right],
\eeq
where $$\varpi_{p}=\int_{0}^{\frac{1}{2}}\frac{[(1+\rho)^{\gamma^{\ast}}-(1-\rho)^{\gamma^{\ast}}]^{p}}{\rho}d\rho$$  and $$C_{2}=\max_{r\in[1/2,1]}\frac{[(1+\rho)^{\gamma^{\ast}}-(1-\rho)^{\gamma^{\ast}}]^{p}}{\rho}.$$
Since

$$
\begin{cases}
\displaystyle\lim_{\rho\rightarrow0^{+}}\frac{[(1+\rho)^{\gamma^{\ast}}-(1-\rho)^{\gamma^{\ast}}]^{p}}{\rho}=0
& \mbox{if } p\in(1,2),\\
\displaystyle\lim_{\rho\rightarrow0^{+}}\frac{(1+\rho)^{\gamma^{\ast}}-(1-\rho)^{\gamma^{\ast}}}{\rho}=2\gamma^{\ast} &\mbox{if } p=1,\\
\displaystyle
\lim_{\rho\rightarrow0^{+}}\left\{\frac{[(1+\rho)^{\gamma^{\ast}}-(1-\rho)^{\gamma^{\ast}}]^{p}}{\rho}\rho^{1-p}\right\}
=(2\gamma^{\ast})^{p}, & \mbox{if } p\in(1/\gamma^{\ast},1),
\end{cases}
$$
we see that $\varpi_{p}$ is convergent, which, together with (\ref{eq-2.0}), implies that (\ref{eq-th-c1}) holds.

Next, we prove (\ref{eq-th-2}).

By the assumption, $f$ admits the standard decomposition 
$f = h + \overline{g} \in \mathcal{S}_{H}(K)$, 
where 
$h(z) = z + \sum_{n=2}^{\infty} a_{n} z^{n}$ 
and 
$g(z) = \sum_{n=1}^{\infty} b_{n} z^{n}$.
Then, for  $n\in\{2,3,\ldots\}$, we have
\beqq
|a_{n}|&=&\left|\frac{1}{2\pi i}\int_{|\zeta|=r}\frac{h(\zeta)}{\zeta^{n+1}}d\zeta\right|=
\left|\frac{1}{2\pi i}\int_{|\zeta|=r}\frac{f(\zeta)}{\zeta^{n+1}}d\zeta\right|\\
&\leq&r^{-n}M_{1}(r,f),
\eeqq
which, together with (\ref{eq-th-c1}), yields that there is a constant $C_{3}(K)$
such that
\beqq\label{eq-x-11}
|a_{n}|\leq C_{3}(K)\frac{1}{r^{n}(1-r)^{\gamma^{\ast}-1}}.
\eeqq
Therefore,
\beqq
|a_{n}|\leq C_{3}(K)\min_{r\in(0,1)}\left[\frac{1}{r^{n}(1-r)^{\gamma^{\ast}-1}}\right]=C_{3}(K)\left(1+\frac{\gamma^{\ast}-1}{n}\right)^{n}
\left(\frac{n}{\gamma^{\ast}-1}+1\right)^{\gamma^{\ast}-1},
\eeqq
which implies that (\ref{eq-th-2}) holds. In particular, if $K=1$, then ${\mathcal S}={\mathcal S}_{H}(1)$.
It follows from De Branges's result in \cite{B} that $\gamma^{\ast}=2$  and $|a_{n}|\leq n$.
 The proof of this theorem is complete. \qed

\begin{Thm}{\rm (\cite[Theorem 4.4]{Clunie-Small-84})}\label{Thm-G}
Each function $f\in{\mathcal S}_H^{0}$ satisfies the inequality
$$|f(z)|\geq\frac{1}{4}\frac{|z|}{(1+|z|)^{2}},~z\in\mathbb{D}.$$
\end{Thm}

\begin{Thm}{\rm (\cite[Theorem 1]{Small})}\label{Thm-H}
If $f\in{\mathcal S}_H^{0}$, then
$$\frac{1}{2\gamma}\left[1-\left(\frac{1-|z|}{1+|z|}\right)^{\gamma}\right]\leq|f(z)|\leq
\frac{1}{2\gamma}\left[\left(\frac{1+|z|}{1-|z|}\right)^{\gamma}-1\right],~z\in\mathbb{D}.$$
\end{Thm}

\subsection*{The proof of Theorem \ref{thm-02}}
$\mathbf{Case~1}.$ Let $p\in(1/\gamma,2]$.

For $n\in\{1,2,\ldots\}$, let $F_{n}(z)=\left(|f(z)|^{2}+\frac{1}{n}\right)^{\frac{1}{2}}$, $z\in\mathbb{D}$.
By (\ref{eq1.2x}) and (\ref{eq-1g}), we have
\beq\label{eq-1}
r\frac{d}{dr}M_{p}^{p}(r,F_{n})&=&\frac{1}{2}\int_{\mathbb{D}_{r}}\Delta(|F_{n}(z)|^{p})dA(z)\\ \nonumber
&\leq&p\int_{\mathbb{D}_{r}}\left(|f(z)|^{2}+\frac{1}{n}\right)^{\frac{p}{2}-1}\Lambda_{f}^{2}(z)dA(z)\\ \nonumber
&\leq&KpI_{1}(r)+K_{0}pI_{2}(r),
\eeq
where $$I_{1}(r)=\int_{\mathbb{D}_{r}}\left(|f(z)|^{2}+\frac{1}{n}\right)^{\frac{p}{2}-1}J_{f}(z)dA(z)$$
and $$I_{2}(r)=\int_{\mathbb{D}_{r}}\left(|f(z)|^{2}+\frac{1}{n}\right)^{\frac{p}{2}-1}dA(z).$$
Now we come to estimate $I_{1}(r)$ and  $I_{2}(r)$. Elementary calculations give
\beq\label{eq-1-r}
I_{1}(r)&\leq&\int_{\mathbb{D}_{m(r)}}\left(|w|^{2}+\frac{1}{n}\right)^{\frac{p}{2}-1}dA(w)
=\frac{1}{\pi}\int_{0}^{2\pi}\int_{0}^{m(r)}\frac{\rho}{\left(\rho^{2}+\frac{1}{n}\right)^{1-\frac{p}{2}}}d\rho d\theta\\ \nonumber
&=&\frac{2}{p}\left(\left(m^{2}(r)+\frac{1}{n}\right)^{\frac{p}{2}}-\frac{1}{n^{\frac{p}{2}}}\right).
\eeq
On the other hand, by Theorem \ref{Thm-G}, we have

\beq\label{eq-2-r}
I_{2}(r)&=&\int_{\mathbb{D}_{r}}\frac{1}{\left(|f(z)|^{2}+\frac{1}{n}\right)^{1-\frac{p}{2}}}dA(z)
\leq\int_{\mathbb{D}_{r}}\frac{dA(z)}{\left[\frac{1}{n}+\frac{1}{16}\frac{|z|^{2}}{(1+|z|)^{4}}\right]^{1-\frac{p}{2}}}\\ \nonumber
&\leq&2^{8-4p}\int_{\mathbb{D}_{r}}\frac{1}{|z|^{2-p}}dA(z)=2^{9-4p}\int_{0}^{r}\frac{1}{\rho^{1-p}}d\rho\\ \nonumber
&=&\frac{2^{9-4p}}{p}r^{p}.
\eeq
It follows from (\ref{eq-1}), (\ref{eq-1-r}), (\ref{eq-2-r}) and   dominated convergence theorem that

\beqq
M_{p}^{p}(r,f)&\leq&2K\int_{0}^{r}\frac{(m(\rho))^{p}}{\rho}d\rho+K_{0}2^{9-4p}\int_{0}^{r}\rho^{p-1}d\rho\\
&=&2K\int_{0}^{r}\frac{(m(\rho))^{p}}{\rho}d\rho+\frac{K_{0}2^{9-4p}}{p}r^{p},
\eeqq
which, together with Theorem \ref{Thm-H} and the similar reasoning as in the proof of (\ref{eq-2.0}), gives
(\ref{eq-th-t1}).

$\mathbf{Case~2}.$ Let $p\in(2,\infty)$.

By (\ref{hh-x9}), we have

\beq\label{hh-x16}
r\frac{d}{dr}M_{p}^{p}(r,f)&\leq&\frac{p^{2}}{2}\int_{\mathbb{D}_{r}}|f(z)|^{p-2}\Lambda_{f}^{2}(z)dA(z)\\ \nonumber
&\leq&\frac{p^{2}K}{2}\int_{\mathbb{D}_{r}}|f(z)|^{p-2}J_{f}(z)dA(z)+
\frac{p^{2}K_{0}}{2}\int_{\mathbb{D}_{r}}|f(z)|^{p-2}dA(z)\\ \nonumber
&=&\frac{p^{2}K}{2}\int_{\mathbb{D}_{m(r)}}|w|^{p-2}dA(w)+\frac{p^{2}K_{0}}{2}\int_{\mathbb{D}_{r}}|f(z)|^{p-2}dA(z),
\eeq
which, together with Theorem \ref{Thm-H}, gives that

\beq\label{eq-dt-0}
M_{p}^{p}(r,f)
\leq pKI_{3}(r)+
\frac{p^{2}K_{0}}{2\gamma}I_{4}(r),
\eeq
where $$I_{3}(r)=\int_{0}^{r}\frac{(m(\rho))^{p}}{\rho}d\rho$$
and
$$I_{4}(r)=\frac{1}{r}\int_{0}^{r}\frac{[(1+\rho)^{\gamma}-(1-\rho)^{\gamma}]^{p-2}}{(1-\rho)^{\gamma(p-2)}}d\rho.$$

We first estimate $I_{3}(r)$. By Theorem \ref{Thm-H}, we have

\beq\label{eq-dt-1}
I_{3}(r)&\leq& \frac{1}{2\gamma}\int_{0}^{r}\frac{\left[\left(\frac{1+\rho}{1-\rho}\right)^{\gamma}-1\right]^{p}}{\rho}d\rho
\leq\frac{C_{1}}{2\gamma}
\int_{0}^{r}\frac{1}{(1-\rho)^{\gamma p}}d\rho\\ \nonumber
&=&\frac{C_{1}}{2\gamma(\gamma p-1)}\left(\frac{1}{(1-r)^{\gamma p-1}}-1\right),
\eeq
where $C_{1}$ is defined in (\ref{tt-1}).

Next, we estimate $I_{4}$.
Since $ I_{4}'(r) \geq 0 $ for all $ r \in (0,1) $, it follows that $ I_{4}(r) $ is monotonically increasing on the interval $ (0,1) $
and $$\lim_{r\rightarrow0^{+}}I_{4}(r)=0.$$
Without loss of generality, we assume $r\in[1/2,1)$. 
For $\rho\in[0,r]$, let $$\tau(\rho)=\frac{[(1+\rho)^{\gamma}-(1-\rho)^{\gamma}]^{p-2}}{(1-\rho)^{\gamma(p-2)}}.$$
Then  $\tau(\rho)$ is monotonically increasing with respect to  $\rho$  on the interval $[0, r]$.
Consequently,
\beq\label{eq-dt-2}
I_{4}(r)&=&\frac{1}{r}\int_{0}^{\frac{1}{2}}\tau(\rho)d\rho+
\frac{1}{r}\int_{\frac{1}{2}}^{r}\tau(\rho)d\rho\\ \nonumber
&\leq&2\int_{0}^{\frac{1}{2}}\tau\left(\frac{1}{2}\right)d\rho+
2\int_{\frac{1}{2}}^{r}\frac{[(1+1)^{\gamma}-(1-1)^{\gamma}]^{p-2}}{(1-\rho)^{\gamma(p-2)}}d\rho\\ \nonumber
&=&2^{\gamma(p-3)}\left(3^{\gamma}-1\right)+2^{1+\gamma(p-2)}I_{5}(r),
\eeq
where %$$I_{5}(r)=\int_{\frac{1}{2}}^{r}\frac{d\rho}{(1-\rho)^{\gamma(p-2)}}.$$

%we see that
%\be\label{eq-dt-2}
%I_{4}\leq C_{4}\int_{0}^{r}\frac{d\rho}{(1-\rho)^{\gamma(p-2)}}\leq C_{4}\frac{1}{(1-r)^{\gamma(p-2)}},
%\ee

\beqq
I_{5}(r)&=&\int_{\frac{1}{2}}^{r}\frac{d\rho}{(1-\rho)^{\gamma(p-2)}}\\
&=&
\begin{cases}
\displaystyle \frac{1}{\gamma(p-2)-1}\left(\frac{1}{(1-r)^{\gamma(p-2)-1}}-2^{\gamma(p-2)-1}\right)
& \mbox{if } \gamma(p-2)\neq1,\\
\displaystyle \log\frac{1}{1-r}-\log2 &\mbox{if } \gamma(p-2)=1.
\end{cases}
\eeqq
Combining (\ref{eq-dt-0}), (\ref{eq-dt-1}) and (\ref{eq-dt-2}) gives  (\ref{eq-th-t1}).

%\beqq
%M_{p}^{p}(r,f)&\leq&pK\int_{0}^{r}\frac{(m(\rho))^{p}}{\rho}d\rho+
%\frac{p^{2}K_{0}}{2r\gamma}\int_{0}^{r}\frac{[(1+\rho)^{\gamma}-(1-\rho)^{\gamma}]^{p-2}}{(1-\rho)^{\gamma(p-2)}}d\rho\\
%&\leq&pK\int_{0}^{r}\frac{(m(\rho))^{p}}{\rho}d\rho+
%\frac{p^{2}K_{0}}{2\gamma}\int_{0}^{r}\frac{[(1+\rho)^{\gamma}-(1-\rho)^{\gamma}]^{p-2}}{\rho}\frac{d\rho}{(1-\rho)^{\gamma(p-2)}}\\
%&\leq&pK\int_{0}^{r}\frac{(m(\rho))^{p}}{\rho}d\rho+
%\frac{p^{2}K_{0}C_{4}}{2\gamma}\int_{0}^{r}\frac{1}{(1-\rho)^{\gamma(p-2)}}d\rho,
%\eeqq
%where $$C_{4}=\sup_{\rho\in(0,1)}\left\{\frac{[(1+\rho)^{\gamma}-(1-\rho)^{\gamma}]^{p-2}}{\rho}\right\}<\infty.$$

By using the similar reasoning as in the proof of (\ref{eq-th-2}), we obtain (\ref{eq-th-t2}). In particular, if $K-1=K_{0}=0$, then $f\in {\mathcal S}_{H}^{0}(1, 0)={\mathcal S}$.
It  follows from De Branges's result in \cite{B} that $\gamma=2$  and $|a_{n}|\leq n$.
 The proof of this
theorem is complete.
\qed\\

{\bf Funding:} This research  was partly supported by the National Science
Foundation of China (grant no. 12571080).\\

{\bf Data Availability Statement:} My manuscript has no associated data.\\

{\bf Conflict of interest:} The author declares that he has no conflict of interest.

%\bigskip
%{\bf Acknowledgements.}\label{acknowledgments}

\bigskip

{\bf Acknowledgements:}  I sincerely thank  two anonymous referees for their thoughtful
comments and constructive suggestions, which have greatly enhanced the quality and readability of this paper. The research of the  author was partly supported by the
National Science Foundation of China (grant no. 12571080), the Key Project of the Natural Science Foundation of Hunan Province (grant no. 2026JJ30002)
and  Guangxi Natural Science Foundation $\#$2026GXNSFFA00640002.

\normalsize

\end{document}